\documentclass[12pt,twoside,reqno]{amsart}
\usepackage{amssymb}
\usepackage{enumerate}
\usepackage[T2A]{fontenc}
\usepackage[utf8]{inputenc}
\usepackage[english]{babel}
\usepackage{url}
\usepackage{color}
\usepackage{hyperref}

\advance\textwidth38mm \advance\hoffset-20mm
\advance\textheight20mm \advance\voffset-11mm
\allowdisplaybreaks[2]
\theoremstyle{plain}
\newtheorem{thm}{Theorem}[section]
\newtheorem{cor}[thm]{Corollary}
\newtheorem{lem}[thm]{Lemma}
\newtheorem{proposition}[thm]{Proposition}
\theoremstyle{remark}
\newtheorem{rem}[thm]{Remark}
\newtheorem{defn}[thm]{Definition}
\newtheorem{exa}[thm]{Example}
\theoremstyle{plain}
\numberwithin{equation}{section}
\newcommand\sign{\operatorname{sign}}
\newcommand\diag{\operatorname{diag}}
\newcommand\Z{{\scriptscriptstyle\mathbb{Z}}}
\newcommand\Zinfty{\infty}
\newcommand\bigcdot{{\Large\boldmath$\cdot$}}
\newcommand\doi[1]{DOI: \url{#1}}
\newcommand\Cdot{{\mskip2mu\cdot\mskip2mu}}
\let\emph=\textit

\begin{document}


\title[Chebyshev systems and Sturm
oscillation theory]{Chebyshev systems and Sturm
oscillation theory for~discrete polynomials}

\author{D.~V.~Gorbachev}
\address{D.\,V.~Gorbachev, Lomonosov Moscow State University, Moscow Сenter of Fundamental and Applied Mathematics, 119991 Moscow, Russia}
\email{dvgmail@mail.ru}

\author{V.~I.~Ivanov}
\address{V.\,I.~Ivanov, Tula State University,
Department of Applied Mathematics and Computer Science,
300012 Tula, Russia;
Lomonosov Moscow State University, Moscow Сenter of Fundamental and Applied Mathematics,
119991 Moscow, Russia}
\email{ivaleryi@mail.ru}

\author{S.~Yu.~Tikhonov}
\address{S.\,Yu.~Tikhonov,
 ICREA, Pg. Llu\'is Companys 23, 08010 Barcelona,
Spain
\\
Centre de Recerca Matem\`atica, Campus de Bellaterra, Edifici C 08193
Bellaterra, Barcelona, Spain,
and Universitat Aut\'onoma de Barcelona}
\email{stikhonov@crm.cat}

\date{\today}



\thanks{The proofs of Theorems~1.4, 1.8, 1.10, 1.14 were supported by
the Russian Science Foundation
(project no.~23-71-30001) at the Lomonosov Moscow State University.
The work of the third author was partially supported by grants PID2023-150984NB-I00, 2021 SGR 00087,
the CERCA Programme of the Generalitat de Catalunya, and by the
Spanish State Research Agency, through the Severo Ochoa and Mar{\'\i}a de
Maeztu Program for Centers and Units of Excellence in R$\&$D
(CEX2020-001084-M).}

\begin{abstract}
We prove an analogue of Chebyshev's alternation theorem for linearly independent discrete functions 
$\Phi_n=\{\varphi_k\}_{k=1}^n$ on the  interval $[0,q]_{\scriptscriptstyle\mathbb{Z}}=[0,q]\cap \mathbb{Z}$. In particular, we establish that the polynomial of best uniform approximation of a discrete function 
admits a Chebyshev alternance set of length $n+1$ if and only if $\Phi_n$ is a Chebyshev $T_{\scriptscriptstyle\mathbb{Z}}$-system.
We also obtain a discrete version of Sturm's oscillation theorem, according to which the number of discrete zeros of the polynomial $\sum_{k=m}^{n}a_k\varphi_k$ is no less than $m-1$ and no more than $n-1$. This implies  that $\Phi_n$ is a $T_{\scriptscriptstyle\mathbb{Z}}$-system and  a discrete Sturm-Hurwitz spectral gap theorem is valid.
As applications, we study the  orthogonal polynomials  with removed largest zeros. We~establish the monotonicity property of coefficients in the Fourier expansions of such polynomials, thereby  strengthening the results of H.~Cohn and A.~Kumar. We apply this to solve a Yudin-type extremal problem for polynomials with spectral gap.
\end{abstract}

\maketitle

\tableofcontents



\section{Introduction}

\subsection{Discrete  Chebyshev system}
Let $C[a,b]$ be the  space of continuous real-valued functions on a finite interval $[a,b]$ equipped with the norm 
$
\|f\|_{\infty}=\max_{t\in [a,b]}|f(t)|.
$
Given linearly independent functions $\{\varphi_k(t)\}_{k=1}^n\subset C[a,b]$, 
 the set of  polynomials is defined by
\[
L_n=\Bigl\{p=\sum_{k=1}^na_k\varphi_k\colon a_k\in\mathbb{R}\Bigr\}.
\]
The best approximation of a continuous function $f$ by 
$L_n$ is given by 
\begin{equation*}
E(f,L_n)_{\infty}=\min_{p\in L_n}\|f-p\|_{\infty}=\|f-p^{*}\|_{\infty}
\end{equation*}
and $p^{*}\in L_n$ is the best approximant. 

The Chebyshev systems play the key role  
in a study of the best approximation of continuous functions by the $n$-dimensional subsets.
Recall that a set of continuous functions $\{\varphi_k\}_{k=1}^n$
is called {a} $T$-system (or {a} Chebyshev system) on  $[a,b]$ if any nontrivial polynomial with respect to this system has at most  $n-1$ distinct zeros on $[a,b]$.
This definition can be equivalently written as follows 
(cf. \cite[Ch.~1, \S\,4]{KS66}): The set  $\{\varphi_k\}_{k=1}^n$ is a $T$-system if and only if, for any sequence  $a\le t_1<\dots<t_n\le b$,
the determinants 
\begin{equation}\label{eq1.2}
\Delta\!\left(\setlength{\arraycolsep}{3pt}
\begin{matrix}
\varphi_1, & \varphi_2, & \dots, & \varphi_n\\
t_1, & t_2, & \dots, & t_n
\end{matrix}\right)
=\begin{vmatrix}
\varphi_1(t_1)&\varphi_2(t_1)&\dots&\varphi_n(t_1)\\
\varphi_1(t_2)&\varphi_2(t_2)&\dots&\varphi_n(t_2)\\
\hdotsfor[2]{4}\\
\varphi_1(t_n)&\varphi_2(t_n)&\dots&\varphi_n(t_n)\\
\end{vmatrix}
\end{equation}
have the same sign.

Let us recall two crucial facts 
 on the best approximations in  $C[a,b]$.

\theoremstyle{plain}\newtheorem*{thmA}{Haar's theorem}
\begin{thmA}[{\cite[Th. 3.4.6]{La72}}]
The best approximant $p^{*}\in L_n$ of a function  $f\in C[a,b]$ is unique if and only if $\{\varphi_k\}_{k=1}^{n}\subset C[a,b]$ is a Chebyshev system.
\end{thmA}

\theoremstyle{plain}\newtheorem*{thmB}{Chebyshev's theorem}
\begin{thmB}[{\cite[Th. 3.4.7]{La72}}]
Let $\{\varphi_k\}_{k=1}^{n}\subset C[a,b]$ be the 
Chebyshev system. A polynomial  $p^{*}\in L_n$ is the best approximant of a function  $f\in C[a,b]$ if and only if there exists an ordered set of  $n+1$ points  $a\le t_1<\dots<t_{n+1}\le b$ such that 
\begin{equation}\label{eq1.3}
\begin{aligned}
\textup{(i)}&\quad |p^{*}(t_i)-f(t_i)|=\|p^{*}-f\|_{\infty},\quad i=1,\dots,n+1,\qquad \\
\textup{(ii)}&\quad p^{*}(t_i)-f(t_i)=-(p^{*}(t_{i+1})-f(t_{i+1})),\quad i=1,\dots,n.
\end{aligned}
\end{equation}
\end{thmB}
Recall that an ordered set of $\{t_i\}$ satisfying properties \eqref{eq1.3} is called the Chebyshev alternance set of length $n+1$.

We are interested in analogues of Haar's
and Chebyshev's theorems for the best uniform approximation of discrete functions defined on 
  ordered sets of points, rather than on an interval. 
  Note that some basic results 
  can be partially derived
from well-known results on the best uniform approximation of continuous functions on an arbitrary compact set (see, e.g.,   \cite{Du79}, \cite[Ch.~1, \S\S\,2, 5]{DS08},~\cite[Ch.~II]{GK02}, \cite[Ch.~3]{La72}).


For $m,n\in\mathbb{Z}$, $m\le n$,
we define
$[m,n]_{\Z}=\{m,m+1,\dots,n\}$.
 For  $q\in\mathbb{N}$, let\footnote{Here, similarly to the continuous case,
we use the notation $C[0,q]_{\Z}$.} 
$C[0,q]_{\Z}=\{f\colon [0,q]_{\Z}\to \mathbb{R}\}$ be the linear space (of the dimension $q+1$) of real-valued discrete functions on $[0,q]_{\Z}$ equipped with the norm 
$
\|f\|_{\Zinfty}=\max_{\nu\in [0,q]_{\Z}}|f(\nu)|.
$ 
The point  $\nu$ is called \textit{a zero} of the discrete function  $f$ if\footnote{Following 
P.~Hartman \cite{Ha78}
such zero is sometimes called generalized; see, e.g.,
\cite[Definition~1.2.10]{ABGR05} or
\cite[Definition~7.8]{El05}.
} either
\[
f(\nu)=0\ \ 
\text{for}\ \ \nu \in [0,q]_{\Z}
\quad \text{or}\quad f(\nu-1)f(\nu)<0\ \ \text{for}\ \ \nu \in [1,q]_{\Z}.
\]
It will be important to distinguish these zeros and  we will call them the zeros of the \textit{first} or \textit{second} type, respectively.
Denote by $N(f,[m,n]_{\Z})$ the number of zeros (of both types) of $f$ on $[m,n]_{\Z}$.
Set also $N(f)=N(f,[0,q]_{\Z})$. Moreover, let $N_0(f)$ be the number of zeros of the first type of $f$ on 
$[0,q]_{\Z}$. Clearly, $N_{0}(f)\le N(f)$.

F.~Gantmacher and M.~Krein~\cite[Ch.~II]{GK02} suggested the following characteristics
of oscillatory properties of discrete functions. 
Let $S^{-}(f)$ (respectively, $S^{+}(f)$) be the least (the largest) number of sign changes\footnote{As usual, $f$ has a sign change on $[0,q]_{\Z}$ if $f(\nu-1)f(\nu)<0$ for some $\nu\in [1,q]_{\Z}$.} of $f$ on $[0,q]_{\Z}$ after replacing all zero values of $f$ by arbitrary 
nonzero values. 
We claim that 
\[
S^{-}(f)\le N(f)\le S^{+}(f),
\]
with the discussion and proof given in Section~\ref{sec3} (see Lemma~\ref{lem7}).

We study the problem of
the characterization and uniqueness of
the best uniform approximants of discrete functions. Let 
$\{\varphi_k\}_{k=1}^{n}\subset C[0,q]_{\Z}$ be 
linearly independent discrete functions and 
$L_n$
 be its linear span over
$\mathbb{R}$.
As in the continuous case,
\[
E(f,L_n)_{\Zinfty}=\min_{p\in L_n}\|f-p\|_{\Zinfty}=\|f-p^{*}\|_{\Zinfty}
\]
is the best uniform approximation of $f$ by $L_n$ and $p^{*}\in L_n$ is the best approximant.


Since, for discrete functions one can distinguish between two types of zeros, we define two corresponding 
 discrete Chebyshev's systems.
Let $n\in [1,q+1]_{\Z}$. The set of discrete functions $\{\varphi_k\}_{k=1}^{n}\subset C[0,q]_{\Z}$ is called the 
$T_{\Z}$-system (respectively, $T_{0}$-system) if for any nontrivial polynomial 
\begin{equation}\label{eq1.4}
p(\nu)=\sum_{k=1}^{n}a_k\varphi_k(\nu)
\end{equation}
of degree at most $n$ with real coefficients we have $N(p)\le n-1$ ($N_0(p)\le n-1$). 
 It is clear that any $T_{\Z}$-system is also $T_0$-system but the inverse statement is not valid, see the example in Section~\ref{sec2}.

Let $p^{*}\in L_n$ be the best uniform approximant of the function $f\in C[0,q]_{\Z}$ and set 
\[
e(\nu)=p^{*}(\nu)-f(\nu), \quad \nu \in [0,q]_{\Z}.
\]
Due to compactness of  $[0,q]_{\Z}\subset \mathbb{R}$, 
 we observe that  an analogue of Haar's theorem is valid only for 
 $T_0$-systems.

\begin{thm}[{\cite[Th.~3.4.6]{La72}}]
The  best uniform approximant 
$p^*\in L_n$
of the function $f\in C[0,q]_{\Z}$ is unique if and only if the set $\{\varphi_k\}_{k=1}^{n}\subset C[0,q]_{\Z}$ is a 
$T_0$-system.
\end{thm}


 We note that a (nontrivial) polynomial \eqref{eq1.4} vanishes at $n$ integer points $0\le\nu_1<\dots<\nu_n\le q$ if and only if the determinant \begin{equation}\label{eq1.5}
\Delta\!\left(\setlength{\arraycolsep}{3pt}
\begin{matrix}
\varphi_1, & \varphi_2, & \dots, & \varphi_n\\
\nu_1, & \nu_2, & \dots, & \nu_n
\end{matrix}\right)
\end{equation}
defined by 
\eqref{eq1.2}
is zero.
Therefore, we arrive at the following obvious statement, 
which provides a criterion for $T_0$-systems
(see, e.g.,  \cite[Ch.~VII, \S\,1]{KS66}).

\begin{proposition}\label{lem1}
The set $\{\varphi_k\}_{k=1}^{n}\subset C[0,q]_{\Z}$ is a 
$T_0$-system if and only if the determinants~\eqref{eq1.5} are nonzero
for any 
 integers $0\le\nu_1<\dots<\nu_n\le q$.
\end{proposition}

A straightforward  analogue  of  Chebyshev's theorem for 
$T_0$-systems
reads as follows.
\begin{thm}[{\cite[Th.~3.4.7]{La72}}]\label{theo2}
Let the set $\{\varphi_k\}_{k=1}^{n}\subset C[0,q]_{\Z}$ be a 
$T_0$-system. The polynomial $p^{*}\in L_n$ is the best uniform approximant of a function $f\in C[0,q]_{\Z}$ if and only if there exist distinct
points $\nu_1, \dots,\nu_{n+1}\in [0,q]_{\Z}$,
$\varepsilon_1,\ldots,\varepsilon_{n+1}\in \{-1,1\}$, and positive numbers $\rho_1,\ldots,\rho_{n+1}$, satisfying the conditions
\begin{equation*}
\begin{aligned}
\textup{(i)}&\quad \varepsilon_i e(\nu_i)=\|e\|_{\Zinfty},\quad i=1,\dots,n+1,\quad e=p^{*}-f,\\
\textup{(ii)}&\quad \sum_{i=1}^{n+1}\varepsilon_i\rho_i p(\nu_i)=0\quad 
\text{for all}\quad p\in L_n.
\end{aligned}
\end{equation*}
\end{thm}

We note that the set of points $\nu_1,\dots,\nu_{n+1}$ 
in Theorem~\ref{theo2} 
may not form the Chebyshev alternance set
as we will see in Example \ref{example2.7}.

Our first main result provides a characterization of the discrete set $\{\varphi_k\}_{k=1}^{n}\subset C[0,q]_{\Z}$ such that for any $f \in C[0,q]_{\Z}$ the best uniform approximant $p^{*}\in L_n$ has an alternance set of length $n+1$, that is, we obtain
a complete analogue of Chebyshev's alternance theorem.



We say that the best uniform approximant $p^{*}\in L_n$ of $f$ admits an alternance set $0\le\nu_1<\dots<\nu_{n+1}\le q$ of length $n+1$ if  the following condition holds (cf.  \eqref{eq1.3}): 
\[
\varepsilon(-1)^{i}e(\nu_i)=\|e\|_{\Zinfty},\quad i=1,\dots,n+1,\quad
\varepsilon=\pm 1.
\]


\begin{thm}\label{theo3} 
Let $n\in [1,q]_{\Z}$. The following  are equivalent$:$
\\
\textup{(a)} 
The set $\{\varphi_k\}_{k=1}^{n}\subset C[0,q]_{\Z}$ is a 
$T_{\Z}$-system.
\\
\textup{(b)} 
For any function $f\in C[0,q]_{\Z}$, the best uniform approximant $p^{*}\in L_n$ of $f$ admits an alternance set of length $n+1$.
\\
\textup{(c)} 
Determinants \eqref{eq1.5} constructed over all sets $0\le\nu_1<\dots<\nu_n\le q$ are nonzero and have the same sign.
\end{thm}
Thus, an analogue of Chebyshev's theorem holds only  for 
$T_{\Z}$-systems.

Using Theorem~\ref{theo3}, 
one can easily construct different
examples of $T_{\Z}$-systems via $T$-systems on an interval.

\begin{cor}\label{cor1}
If $\{\varphi_i(t)\}_{i=1}^n\subset C[0,q]$ is a 
$T$-system on  $[0,q]$, then   
 $\{\varphi_i(\nu)\}_{i=1}^n\subset C[0,q]_{\Z}$ is a 
$T_{\Z}$-system.
\end{cor}

Another important property of 
$T_{\Z}$-systems
 is given below. 

\begin{cor}\label{cor2}

If  $\{\varphi_k\}_{k=1}^n\subset C[0,q]_{\Z}$ is a 
$T_{\Z}$-system, then for any nontrivial polynomial $p \in L_n$ we have 
 $S^{+}(p)\le n-1$.
\end{cor}

{\it Historical comments.}
S.~Karlin and W.J.~Studden \cite[Ch.~VII]{KS66} introduced $T_{\Z}$-systems
via condition (c) and showed that~(c) yields (a). In the matrix form the implication $\text{(b)}\Leftrightarrow \text{(c)}$ can be found in the recent paper \cite{MZO24}.
Corollary \ref{cor1}
 was given in  \cite[Ch.~VII, \S\,1]{KS66} in the particular case when  $T$-systems are generated by  strictly total positive kernels.
Also, Corollary \ref{cor2}
 was proved in \cite[Ch.~VII, \S\,2]{KS66} using property
(c) of Theorem~\ref{theo3}. 
We will give a different proof of Corollary \ref{cor2} based on an oscillation property of discrete polynomials~\eqref{eq2.4}.

We note that the discrete Chebyshev approximation plays an important role in numerical analysis (see, e.g., \cite{Du79}  and \cite{MZO24,ZMT22} for the the best approximation by low-rank matrices). 
In these problems the use of effective algorithms such as the Remez algorithm is crucial. Let us mention that versions of Theorem~\ref{theo2} also allow one to construct effective algorithms for non-Chebyshev systems, see \cite{PK25}.
Moreover, 
C. Dunham
\cite{Du79} 
 developed  a discrete Remez algorithm based on alternance sets.
Thus, Theorem~\ref{theo3} 
completes this study  as it characterizes
{\it all systems} that admit alternation.


For applications, it is important to construct various discrete Chebyshev systems. 
To this end, in the next section we
consider the eigenfunctions of Sturm-Liouville problems.



\subsection{
Sturm oscillation theorem, $T$-systems, and orthogonal polynomials}
We start with the following known fact:
$T$-systems can be obtained using Sturm-Liouville problems. In 1836, C.~Sturm proved the following remarkable result, which was obtained in the same year by J.~Liouville using a different method under less restrictions. See \cite{BH20,Si05,St20} for the historical comments.

\theoremstyle{plain}\newtheorem*{thmC}{Sturm's theorem}
\begin{thmC}
Let $\{V_{k}\}_{k=1}^{\infty}\subset C^{2}[a,b]$ be the set of eigenfunctions associated to eigenvalues $\rho_{1}<\rho_{2}<\dots$ of the following {Sturm-Liouville} problem:
\begin{equation}\label{eq1.7}
\begin{gathered}
(K(t)V'(t))'+(\rho G(t)-L(t))V(t)=0,\quad t\in [a,b],\\ (KV'-hV)(a)=0,\quad (KV'+HV)(b)=0,
\end{gathered}
\end{equation}
where $G,K,L\in C[a,b]$, $K\in C^{1}(a,b)$, $K,G>0$ on $(a,b)$, $h,H\in
[0,\infty]$, and $\rho$ denotes the spectral parameter. Then
\\
\textup{(a)}
Every $k$-th eigenfunction $V_{k}$ has exactly $k-1$
simple zeros in $(a,b)$.
\\
\textup{(b)} For any nontrivial
real polynomial of the form
\begin{equation}\label{eq1.8}
P(t)=\sum_{k=m}^{n}A_{k}V_{k}(t),\quad m,n\in \mathbb{N},\quad m\le n,
\end{equation}
the number of distinct zeros of $P$ on $(a,b)$
is at least $m-1$ and at most
$n-1$. 
\end{thmC}


\begin{rem}\label{rem-SL}
(i) Part (a) is sometimes called the 
weak Sturm oscillation theorem \cite{St20} and proved using comparison theorems (see, e.g., \cite{LS75} and \cite{Si05}).  
Part (b) is a much stronger result
(the 
Sturm oscillation theorem).
As noted by P.~B\'erard and B.~Helffer \cite{BH20}, while part (a) is a well-known statement, part (b) has been almost forgotten. However, for the trigonometric system, it has been known since the early 20th century as the Sturm-Hurwitz theorem:
 \textit{If a real Fourier series has a spectral gap, that is,  $f(t) = \sum_{k=m}^\infty (a_k \cos kt + b_k \sin kt)$, then $f$ has at least $2m$ sign changes over the  period.}
Sturm proved this fact for polynomials (cf. \eqref{eq1.8}), while Hurwitz extended it to the general case.
 Similar results on spectral gap problems are closely related to the Fourier uncertainty principle (see, e.g., 
 \cite{EN04, GIT20,GIT24, Lo77,  MP15, St20, Ul06}) and  signal processing (see, e.g.,~\cite{Ma12}).


(ii) It follows from Sturm's theorem that 
 $\{V_{k}(t)\}_{k=1}^{n}$ 
  form a 
$T$-system on the interval $(a,b)$. By Corollary~\ref{cor1}, for {$a<0$ and $b>q$}, the set $\{V_{k}(\nu)\}_{k=1}^{n}\subset C[0,q]_{\Z}$ is a 
$T_{\Z}$-system.

\end{rem}

Our main goal in this section is 
to
obtain the  Sturm oscillation theorem
for 
  eigenfunctions of a discrete Sturm-Liouville problem on $[0,q]_{\Z}$.
This problem is treated in detail in several monographs, for example,
\cite[Ch.~12]{Ag00}, \cite[Ch.~4]{At64},  \cite[Ch.~7]{KP01}.
 Moreover,  discrete comparison and separation  theorems in various settings 
 have been studied in, e.g., \cite{ABGR05,El05,Si05,Te00}. 
 However, 
an explicit  discrete analogue of   Sturm's theorem has been  an open question.


To formulate the discrete Sturm-Liouville problem, we use infinite tridiagonal Jacobi matrix, which is, in  general, nonsymmetric.
Let $\{\alpha_{l}\}_{l=0}^{\infty}$ be a sequence of real numbers,
$\{\gamma_{l}\}_{l=0}^{\infty}$, $\{\beta_{l}\}_{l=0}^{\infty}$, $\{\rho_l\}_{l=0}^{\infty}$ be
sequences of positive numbers. Define 
\[
J=
\begin{pmatrix}
\alpha_0&\beta_0&0&\dots&0&0&0&\dots\\
\gamma_0&\alpha_1&\beta_1&\dots&0&0&0&\dots\\
\hdotsfor[2]{8}\\
0&0&0&\dots&\gamma_{l-1}&\alpha_{l}&\beta_{l}&\dots\\
\hdotsfor[2]{8}\\
\end{pmatrix}.
\]

Considering the  eigenvalue problem $$JP=\lambda \rho P$$ with
an infinite vector $P=\{P_l(\lambda)\}_{l=0}^{\infty}$ and a diagonal matrix $\rho=\diag(\{\rho_l\}_{l=0}^{\infty})$, we arrive at the  Sturm-Liouville problem in the
 recurrence relation
form:
\begin{equation}\label{eq1.10}
\begin{gathered}
\gamma_{l-1}P_{l-1}(\lambda)+\alpha_lP_l(\lambda)+\beta_{l}P_{l+1}(\lambda)=\lambda
\rho_{l}P_l(\lambda),\quad l\in \mathbb{Z}_{+},\quad \lambda\in\mathbb{C},\\
P_{-1}(\lambda)=0,\quad P_0(\lambda)=1,
\end{gathered}
\end{equation}
where $\gamma_{-1}>0$ is arbitrary. It is clear that the solution of this problem is
given by a family of algebraic polynomials
{$P_l(\lambda)$ of degree $l$} with positive leading coefficients. 



For a non-decreasing function $\mu(\lambda)$ on $\mathbb{R}$, 
the set
\[
\mathcal{S}(\mu)=\{\lambda\colon \mu(\lambda+\delta)-\mu(\lambda-\delta)> 0\ \text{for all}\ \delta>0\}
\]
is called the spectrum of $\mu$. 
To attack Problem \eqref{eq1.10}, we apply  the Favard theorem (see
\cite[Th.~II. 6.4]{Ch78}), which leads to the following statement.

\begin{thm}\label{theo4}
\textup{(a)} Polynomials $P_l(\lambda)$   
are orthogonal with respect to
a positive measure $\mu$ on~$\mathbb{R}$, defined by a non-decreasing function of bounded variation with infinite spectrum, and
\begin{equation}\label{eq1.12}
d_{l}\int_{\mathbb{R}}P_l(\lambda)P_m(\lambda)\,d\mu(\lambda)=\delta_{lm},\quad
d_0=\rho_0,\quad d_l=\rho_{l}\,\frac{\beta_0\cdots \beta_{l-1}}{\gamma_0\cdots
\gamma_{l-1}},\quad l\ge 1.
\end{equation}

\textup{(b)} A measure $\mu$ has a finite support $[a,b]$, the endpoints of which will be the limit points of zeros of the orthogonal polynomials $P_l(\lambda)$, if and only if the sequences
\begin{equation}\label{a-g}
\frac{\alpha_{l}}{\rho_{l}},\quad
\frac{\gamma_{l-1}\beta_{l-1}}{\rho_{l-1}\rho_{l}}
\end{equation}
are 
 bounded.
\end{thm}

Note that positivity of coefficients $\gamma_l$, $\beta_l$, $\rho_l$ is a natural condition in Theorem \ref{theo4}, see Theorem~II.6.4 in \cite{Ch78}.



Throughout the paper we
 assume that
$
q\in \mathbb{Z}_+$ {and} $\eta\in\mathbb{R}.$
Consider the polynomial
\begin{equation}\label{eq1.13}
\widetilde{P}_{q+1}(\lambda)=P_{q+1}(\lambda)-\eta P_{q}(\lambda)
\end{equation}
and define its zeros by
$
\lambda_{q+1}<\dots<\lambda_{1},
$
see
\cite[Ch.~III, \S\,3.3]{Sz59}.
We also  define the following  set of discrete functions: 
\begin{equation}\label{eq1.14}
\psi_k(\nu)=P_{\nu}(\lambda_k),\quad k\in [1,q+1]_{\Z},\quad \nu\in [0,q]_{\Z}.
\end{equation}
By virtue of \eqref{eq1.10}, $\lambda_k$ and $\psi_k(\nu)$ are eigenvalues and eigenfunctions of the discrete Sturm-Liouville problem
\begin{equation}\label{eq1.15}
\begin{gathered}
\gamma_{\nu-1}\psi(\nu-1)+\alpha_{\nu}\psi(\nu)+\beta_{\nu}\psi(\nu+1)=\lambda
\rho_{\nu}\psi(\nu),\quad \nu\in [0,q]_{\Z},\\
\psi(-1)=0,\quad \psi(0)=1,\quad \psi(q+1)=\eta\psi(q).
\end{gathered}
\end{equation}
Introducing the Jacobi matrix
\begin{equation*}
J_{q+1}=
\begin{pmatrix}
\alpha_0&\beta_0&0&\dots&0&0&0&0&\\
\gamma_0&\alpha_1&\beta_1&\dots&0&0&0&0&\\
\hdotsfor[2]{8}\\
0&0&0&\dots&0&\gamma_{q-2}&\alpha_{q-1}&\beta_{q-1}\\
0&0&0&\dots&0&0&\gamma_{q-1}&\alpha_{q}+\eta \beta_{q}\\
\end{pmatrix},
\end{equation*}
 problem \eqref{eq1.15} corresponds to the  eigenvalue problem $$J_{q+1}\psi=\lambda \rho \psi.$$ 

Sturm's theory on the zeros of discrete polynomials {with respect to the set} $\{\psi_k(\nu)\}_{k=1}^{q+1}$ (see \eqref{eq1.14}) 
is similar to that in the continuous case
 (see \eqref{eq1.7}), 
with the additional use of the quantities
$S^{-}(f)$ and $S^{+}(f)$.

Before 
stating  an analogue of item (a) of  Sturm's theorem given in Theorem \ref{theo5} below,
 we mention the following result  by B. Simon.\footnote{A similar result can be found in
\cite[Theorem~4.3.5]{At64}.}
\begin{thm}[{\cite[Th.~2.3]{Si05}}]
For $n\in [1,q]_\Z,$ there holds
\[
|\{j\in [1,n]_\Z\colon \lambda<\lambda_{j,n}\}|=
|\{\nu\in [1,n]_\Z\colon \sign P_{\nu-1}(\lambda)\ne 
\sign P_{\nu}(\lambda)\}|,
\]
where 
$\lambda_{l,l}<\dots<\lambda_{1,l}$ are the zeros of the polynomials $P_l(\lambda)$,
$\lambda\ne \lambda_{j,l}$ for $1\le j\le l\le n$,
and $|I|$ denotes the cardinality of a finite set $I$.
\end{thm}

Taking into account that the function $\psi_k$ cannot have two consecutive zeros of the first type,
this theorem implies the following result.
\begin{thm}\label{theo5}
For all $k\in [1,q+1]_{\Z},$
\[
S^{-}(\psi_k)=S^{+}(\psi_k)=N(\psi_k)=k-1.
\]
\end{thm}

Note that Theorem~\ref{theo5} also follows from 
Theorem~1 in \cite[Ch.~II]{GK02} whose  proof 
 based on the use of the concepts of lines and nodes 
for the piecewise linear function 
associated with Jacobi matrices (see 
Definition~2 in \cite[Ch.~II]{GK02}). 
Moreover, the proof of the fact 
that $N(\psi_k)=k-1$ 
 based on the discrete Sturm separation theorem
can be found in \cite[Theorem~7.6]{KP01}.
The corresponding results  for self-adjoint systems are given in 
\cite[Theorem~2]{DK07},
\cite[Corollary~5.79]{DEH19}, \cite[Theorem~1]{Kr02}.
In Section~\ref{sec3}, we present a simple proof of Theorem~\ref{theo5} based on the property of the interlacing of zeros of  orthogonal polynomials (see \eqref{eq3.1}). 

Our main result in this section is the discrete counterpart of  item (b) of  Sturm's theorem, cf. 
Remark~\ref{rem-SL}.

\begin{thm}\label{theo6} 
Let  integers $m$ and $n$  satisfy $1\le m\le n\le q+1$. For any nontrivial polynomial
\begin{equation}\label{eq1.17}
 V(\nu)=\sum_{k=m}^{n}a_k\psi_k(\nu),
\end{equation}
we have
\begin{equation}\label{eq1.17+}
m-1\le S^{-}(V)\le N(V)\le S^{+}(V)\le n-1.
\end{equation}
\end{thm}

The proof of Theorem~\ref{theo6} is based on Theorem~\ref{theo5} and a generalization of  Liouville's method (see \cite{BH20}).
Note also that the inequalities $m-1\le S^{-}(V)$ and $S^{+}(V)\le n-1$ for $\eta=0$ and $\rho_l=1$ were proven in \cite[Ch.~II, Th.~6]{GK02} under the assumption that the Jacobi
matrix $J_{q+1}$  is oscillatory. 
Recall that a square matrix is called oscillatory if it is totally non-negative, and some power of it is totally positive (see \cite[Ch.~II, Definition~4]{GK02}). In our case, in general, the matrix $J_{q+1}$ is not oscillatory and it becomes oscillatory only if all its principal minors are positive \cite[Ch.~II, Th.~11]{GK02}.
It is worth mentioning 
 that the estimate 
$N(V)\le  n-1$ was 
mentioned  in \cite[
Problem~2,
p.~539]{At64}.

\begin{cor}
For any $n\in [1,q+1]_{\Z}$, the set $\{\psi_k(\nu)\}_{k=1}^{n}$ is a  
$T_{\Z}$-system.
\end{cor}

To conclude, we state  the discrete spectral gap theorem (cf. Remark \ref{rem-SL}).
\begin{cor}
Let $m\in [1,q+1]_{\Z}$ and  $f\in C[0,q]_{\Z}$ be such that 
$f(\nu)=\sum_{k=m}^{q+1}a_k\psi_k(\nu)$.
Then $m-1\le S^{-}(f)\le N(f)$.
\end{cor}


\subsection{Monotonicity property of coefficients of polynomials with removed largest  zeros}
Consider  polynomial \eqref{eq1.13} and divide it into $\prod_{j=1}^{m+1} (\lambda-\lambda_{j})$ with    $m+1$ largest zeros $\lambda_{j}$, $m\le q$. 
Taking into account \eqref{eq1.12}, we expand it
 into the Fourier sum
\begin{equation}\label{eq1.18}
\frac{\widetilde{P}_{q+1}(\lambda)}{(\lambda-\lambda_{1})\cdots (\lambda-\lambda_{m+1})}=
\sum_{l=0}^{q-m}d_{l}a_{l}P_{l}(\lambda).
\end{equation}

The question arises whether all coefficients $a_{l}\ge 0$. It turns out that the answer is affirmative. The case $m=0$ is simple and corresponds to the Christoffel-Darboux formula, while the case
$m=1$ was proven in \cite{GI00}. In full generality, this non-trivial fact was proven by H.~Cohn and A.~Kumar \cite[Th. 3.1]{CK07} 
while working on  the problem on discrete energy on the {Euclidean} sphere. 
 
 In this section, using the 
 discrete Sturm theorem,
 we take one step further and 
 establish the  monotonicity property of the Fourier coefficients for all $m$, which also implies that $a_{l}\ge 0$. 
This fact is  crucial to deal
with a spectral gap problem in Section \ref{subsec-sgp}.

Also, we note that the related  
 polynomial 
  $p_m(\lambda)=\frac{P_{q+1}^2(\lambda)}{(\lambda-\lambda_{1})\cdots (\lambda-\lambda_{m+1})}$  appeared in various extremal problems.
 For example, S. Bernstein used $p_0$ to disprove the existence of Chebyshev's quadrature formula.
 Moreover, 
 this polynomial is an
 extremizer
 in the following question (see  \cite{Ba84,MU84}):  find
 $$\inf_{a_k, b_k} \operatorname{mes}\,\Bigl\{x\in (-\pi,\pi]\colon \sum_{k=1}^n (a_k\cos kx+b_k\sin kx)>0\Bigr\}.$$
  V.~Yudin used $p_m$  in problems of multiple covering of the torus \cite{Yu02} and the Euclidean sphere~\cite{Yu05}.

  Our main result here is the following theorem.
 




\begin{thm}\label{theo7}
Let $q\in\mathbb{N}$ and $m\in [0,q]_{\Z}$. Let also $(a,b)$ be an  interval containing all the zeros of the polynomial $P_{q }(\lambda)$ and $\eta_{b}=P_{q+1}(b)/P_{q}(b)$. 

The following inequalities hold for the coefficients $a_l$ in the expansion \eqref{eq1.18}$:$
\\
\textup{(a)} if either $m\ge1$ and $\eta\in\mathbb{R}$  or $m=0$ and $\eta<\eta_{b}$, then
\begin{equation}\label{eq1.19}
\frac{a_0}{P_0(b)}>\frac{a_1}{P_1(b)}>\cdots>\frac{a_{q-m}}{P_{q-m}(b)}>0;
\end{equation}
\textup{(b)} if $m=0$ and $\eta=\eta_{b}$, then
\begin{equation}\label{eq1.20}
\frac{a_0}{P_0(b)}=\frac{a_1}{P_1(b)}=\cdots=\frac{a_{q}}{P_{q}(b)}>0;
\end{equation}
\textup{(c)} if $m=0$ and $\eta>\eta_b$, then
\begin{equation}\label{eq1.21}
0<\frac{a_0}{P_0(b)}<\frac{a_1}{P_1(b)}<\cdots<\frac{a_{q}}{P_{q}(b)}.
\end{equation}
\end{thm}



For polynomials normalized by $P_{l}(1)=1$,
this result 
 has a simpler form.
 
\begin{cor}\label{cor-theo7}
Let  $P_{l}(\lambda)$ be orthogonal  on $[-1,1]$ and $P_{l}(1)=1$.
 Then 
\[
{\begin{cases}
a_0>a_1>\cdots>a_{q-m}>0,&\text{if either $m=0$ and $\eta<1$ or $m\ge1$ and $\eta\in\mathbb{R}$},
\\
a_0=a_1=\cdots=a_{q}>0,&\text{if $m=0$ and $\eta=1$},\\
0<a_0<a_1<\cdots<a_{q},&\text{if $m=0$ and $\eta>1$}.
\end{cases}}
\]
\end{cor}

\subsection{An extremal problem for polynomials with spectral gap}\label{subsec-sgp}
{Let $m,n\in\mathbb{Z}_+$, $m<n$.}
Let also $\mu$ be the measure on $[-1,1]$ from Theorem~\ref{theo4} and $\{U_{l}\}$ be the set of polynomials orthogonal on $[-1,1]$ with respect to $d\mu$ with normalization $U_{l}(1 )=1$,
$d_l\int_{-1}^1U_lU_m\,d\mu=\delta_{lm}$.

We consider  the following problem: find 
$\lambda\in [-1,1]$ and
{an algebraic polynomial of degree~$n$} on $[-1,1]$ with non-negative Fourier coefficients 
in $\{U_{l}\}$
{having a spectral gap of the length~$m$ (equivalently, with the first $m$ moments being 
zero)}, which preserves its sign over the largest possible interval $[-1,\lambda]$ .

Let $\Pi$ be the set of real algebraic polynomials and $\Pi_n$ consist of those  
 of degree at most $n$.
By $\Pi_{+}(\{U_{l}\})$ we denote
a subset of  $p\in \Pi$ such that 
all coefficients $a_{l}$
in the Fourier expansion 
\begin{equation*}
p(t)=\sum_{l}d_{l}a_{l}U_{l}(t),\quad
a_{l}=a_l(p)=\int_{-1}^{1}p(t)U_{l}(t)\,d\mu(t), 
\end{equation*}
are  non-negative.


Set 
\begin{equation*}
K_{n}(\mu,m)=\Big\{p\in\Pi_n\cap\Pi_{+}(\{U_l\})\colon
\mu_i(p)=0,\ i=0,\dots,m-1,\ \mu_{m}(p)\ge 0\Big\},
\end{equation*}
where $\mu_i(p)=\int_{-1}^{1}t^ip(t)\,d\mu(t)$ is the $i$th moment of the polynomial $p$. 
Note that 
the conditions $\mu_i(p)=0$ are equivalent to the fact that $a_i(p)=0$, $ i=0,\dots,m-1.$

Our goal
is
to find the quantity
\[
B_n(\mu,m)=\sup_{p\in K_{n}(\mu,m)\setminus \{0\}}\{\lambda\in
[-1,1]\colon (-1)^{m-1}p(t)\ge 0,\ -1\le t\le\lambda\}.
\]
This question turns out to play a crucial role in the study of spherical codes and designs.
For $m=1$ it was posed by V.~Yudin \cite{Yu97}. 
Specifically, for $m=1$ and $n=2q$, he 
found extremizers for this problem, though the positivity of the Fourier coefficients was not proven.
Further results in this direction can be found in \cite{GI00,GIT20,GIT24}; see also 
\cite{CK07, Le83}.



Let $\{U_l^{(1)}(t)\}_{l=0}^{\infty}$ be the set of polynomials with normalization $U_l^{(1)}(1) =1$, orthogonal on $[-1,1]$ with respect to the measure {$d\mu^{(1)}(t)=(1+t)\,d\mu(t)$}. The zeros of the polynomials $U_{q+1}(t)$ and $U_{q+1}^{(1)}(t)$, numbered in the descending order, are denoted by $t_k$ and $t_k^{(1)}$, $k=1,\dots,q+1$, respectively. 

{We say  that the set of orthogonal polynomials $\{V_l\}$ satisfies \textit{the Krein property} if $V_mV_n\in\Pi_{+}(\{V_l\})$, that is, $V_{m}V_{n}=\sum_{k=|m-n|}^{m+n}c_{m,n,k}V_{k}$ with $c_{m,n,k}\ge 0$,}
for any $m,n\in\mathbb{Z}_+$, see \cite{Le83}.

\begin{thm}\label{theo8}
Let $m\le q$ and the set $\{U_l\}$ satisfy the Krein property. 
\\
\textup{(a)}
If $n=2q-m+1$, then
\begin{equation}\label{eq1.24}
B_n(\mu,m)=t_{m+1}
\end{equation}
and the unique extremal polynomial  has the form (up to a positive constant)
\begin{equation}\label{eq1.25}
p_{n,1}^{*}(t)=\frac{U_{q+1}^2(t)}{(t-t_{1})\cdots(t-t_{m+1})}.
\end{equation}
\textup{(b)}
If $n=2q-m+2$ and {either} the set $\{U_l^{(1)}\}$ satisfies the Krein property or $\mu(t)$ is an odd function, then
\begin{equation}\label{eq1.26}
B_n(\mu,m)=t_{m+1}^{(1)}
\end{equation}
and the unique extremal polynomial  has the form (up to a positive constant)
\begin{equation}\label{eq1.27}
p_{n,2}^{*}(t)=\frac{(1+t)(U_{q+1}^{(1)}(t))^2}{(t-t_{1}^{(1)})\cdots(t-t_{m+1}^{(1)})}.
\end{equation}
Moreover,  we have 
$\mu_{m}(p_{n,j}^{*})=a_m(p_{n,j}^{*})=0$, $j=1,2$.
\end{thm}

We note that 
a version of Yudin-type problem for a wider class of polynomials $p$ {(without the condition   $p\in\Pi_{+}(\{U_l\})$)} was solved
in \cite{Iv21}. 
Our proof of 
Theorem \ref{theo8}
is based on 
 Corollary~\ref{cor-theo7}.


\subsection{Structure of the paper}

The rest  of the paper is organized as follows. In Section~\ref{sec2}, we provide the proofs of
 Chebyshev's theorem  for 
$T_{\Z}$-systems (Theorem~\ref{theo3}) 
and 
the important
Corollary~\ref{cor2}, which claims that 
 $S^{+}(p)\le n-1$ for a polynomial 
$p \in L_n$, where 
$p(\nu)=\sum_{k=1}^{n}a_k\varphi_k(\nu)$ and 
  $\{\varphi_k\}_{k=1}^n\subset C[0,q]_{\Z}$
  is
a 
$T_{\Z}$-system. 


In Section \ref{sec3}, we 
study the discrete 
Sturm
oscillation theory
and its connection with $T$-systems
and prove Theorems~\ref{theo4}, \ref{theo5}, and \ref{theo6}.

 Section \ref{sec4} is devoted to the proof of Theorem~\ref{theo7}.
In Section \ref{sec5}, 
we study  polynomials with spectral gaps and  solve
Yudin's extremal problem (Theorem~\ref{theo8}).

In the Appendix, we explicitly compute  the determinants \eqref{eq4.1} for  the Jacobi polynomials corresponding to the classical trigonometric systems. 
In particular, we verify that, for all $\nu\in [0,q]_{\Z}$ and $T_{\Z}$-systems, 
$$
\sign 
\Delta\!\left(\setlength{\arraycolsep}{3pt}
\begin{matrix}
\varphi_1, & \dots, &\varphi_{m},& \varphi_{m+1}\\
\nu_{m}, & \dots, &\nu_1, &\nu
\end{matrix}\right)
=\sign 
\Delta\!\left(\setlength{\arraycolsep}{3pt}
\begin{matrix}
\varphi_1, & \dots, &\varphi_{m},& \varphi_{m+1}\\
q, & \dots, &q-m+1, &0
\end{matrix}\right)
\sign \prod_{j=1}^{m}(\nu_j-\nu),
$$
which 
in the general case  will be shown in 
Lemma \ref{lem14}.


\medbreak

\section{Proofs of Theorem~\ref{theo3} and Corollary~\ref{cor2}}\label{sec2}

Recall that 
$q\in \mathbb{Z}_+$ and $\eta\in\mathbb{R}$.
The proof of Theorem~\ref{theo3} is based on 
Lemmas 
\ref{lem2}--\ref{lem6}
below.
First, in addition to the characteristic property of 
$T_0$-systems given by Proposition~\ref{lem1}, we present   another necessary and sufficient condition.

\begin{lem}[{\cite[Proposition~3.4.2]{La72}}]\label{lem2}
A set $\{\varphi_k\}_{k=1}^{n}\subset C[0,q]_{\Z}$ is a 
$T_0$-system if and only if,
for any {integers} $0\le\nu_1<\dots<\nu_n\le q$ and any points $\{y_i\}_{i=1}^n\subset\mathbb{R}$, there is a unique polynomial \eqref{eq1.4}, for which
\[
p(\nu_i)=y_i,\quad i=1,\dots,n.
\]
\end{lem}

Second, we assume that the set $\{\varphi_k\}_{k=1}^{n} \subset C[0,q]_{\Z}$ is a 
$T_0$-system and the determinants~\eqref{eq1.5}
are nonzero
for any integers $0 \leq \nu_1 < \dots < \nu_n \leq q$.

Set, for $0\le\nu_1<\dots<\nu_{n+1}\le q$, 
\begin{align*}
\Delta_{\nu_1}&=\Delta\!\left(\setlength{\arraycolsep}{3pt}
\begin{matrix}
\varphi_1,&\dots,&\varphi_n\\
\nu_2,&\dots,&\nu_{n+1}
\end{matrix}\right),\quad
\Delta_{\nu_{n+1}}=\Delta\!\left(\setlength{\arraycolsep}{3pt}
\begin{matrix}
\varphi_1,&\dots,&\varphi_n\\
\nu_1,&\dots,&\nu_{n}
\end{matrix}\right),
\\
\Delta_{\nu_i}&=\Delta\!\left(\setlength{\arraycolsep}{3pt}
\begin{matrix}
\varphi_1,&\dots,& \varphi_{i-1},& \varphi_{i},&\dots,& \varphi_n\\
\nu_1,&\dots,&\nu_{i-1},&\nu_{i+1},&\dots,& \nu_{n+1}
\end{matrix}\right),
\quad i=2,\dots,n.
\end{align*}

\newcommand\lambdA{\gamma}

\begin{lem}\label{lem3}
Let $0\le\nu_1<\dots<\nu_{n+1}\le q$. 
The linear functional
\begin{equation}\label{eq2.1}
l(p)=\sum_{i=1}^{n+1}\lambdA_{i}p(\nu_i)=0
\end{equation}
for any polynomial $p\in L_n$ if and only if, up to a nonzero constant, $\lambdA_{i}=(-1)^i\Delta_{\nu_i}$, $i=1,\dots,n+1$.
\end{lem}

\begin{proof}
First, we note that  $p\in L_n$ only if  
\[
\Delta\!\left(\setlength{\arraycolsep}{3pt}
\begin{matrix}
p,&\varphi_1, & \varphi_2, & \dots, & \varphi_n\\
\nu_1,&\nu_2, & \nu_3, & \dots, & \nu_{n+1}
\end{matrix}\right)
=0.
\]
Expanding the determinant along the first column, we obtain
\begin{equation}\label{eq2.2}
\sum_{i=1}^{n+1}(-1)^{i}\Delta_{\nu_i}p(\nu_i)=0.
\end{equation}
On the other hand, from \eqref{eq2.1}, for any $a_1,\dots,a_k$,
\[
0=\sum_{i=1}^{n+1}\lambdA_i\sum_{k=1}^{n}a_k\varphi_k(\nu_i)=\sum_{k=1}^{n}a_k\sum_{i=1}^{n+1}\lambdA_i\varphi_k(\nu_i),
\]
and therefore,
\begin{equation}\label{eq2.3}
\sum_{i=1}^{n+1}\lambdA_i\varphi_k(\nu_i)=0,\quad k=1,\dots,n.
\end{equation}
Since  $\Delta_{\nu_i}\neq 0$, all the solutions  $(\lambdA_1,\dots,\lambdA_{n+1})$ 
of the
system \eqref{eq2.3} form a one-dimensional space. Taking into account \eqref{eq2.2}, we arrive at  the statement of the lemma. 
\end{proof}

\begin{lem}\label{lem4}
Let the set $\{\varphi_k\}_{k=1}^{n}\subset C[0,q]_{\Z}$ be a 
$T_{\Z}$-system. The polynomial $p^{*}\in L_n$ is the best uniform approximant of a function $f\in C[0,q]_{\Z}$ if and only if 
$p^{*}$
admits an alternance set of length $n+1$.
\end{lem}

\begin{proof}
The sufficiency  follows from Theorem~\ref{theo2}.
To prove the necessity, according to Lemma~\ref{lem3} and Theorem~\ref{theo2},  we need to show that for a linear functional $l(p)$ on $L_n$ satisfying property \eqref{eq2.1}, there holds $\lambdA_i\lambdA_{i+1}<0$  for $i=1,\dots,n$. We will follow the reasoning given  in the proof of Proposition~3.5.1 in \cite{La72}. Applying Lemma~\ref{lem2}, we construct polynomials $p_i\in L_n$ satisfying 
\[
p_i(\nu_j)=0,\quad j=1,\dots,i-1, i+2,\dots,n+1,\quad p_i(\nu_i)=1.
\]
Since $p_i$ has $n-1$ zeros of the first type, it follows from the definition of the 
$T_{\Z}$-system that {$p_i$} has no zeros of the second type  and, therefore, $p_i( \nu_{i+1})>0$. From this,
\[
l(p_i)=\lambdA_ip_i(\nu_i)+\lambdA_{i+1}p_{i}(\nu_{i+1})=0\quad \text{and}\quad \lambdA_i\lambdA_{i+1}<0,
\]
completing the proof.
\end{proof}

\begin{lem}\label{lem5}
\textup{(a)} For any $f\in C[0,q]_{\Z}$, 
its best uniform approximant
$p^*\in L_n$
 admits an alternance set of length $n+1$ if and only if all determinants \eqref{eq1.5}, constructed over all sequences $0\le\nu_1<\dots<\nu_{n}\le q$, have the same sign.

\textup{(b)} If all determinants \eqref{eq1.5} have the same sign, then any sequence $0\le\nu_1<\dots<\nu_{n+1}\le q$ is an alternance set of length $n+1$ for the polynomial $p^*\in L_n$, which is the best approximant of some function $f\in C[0,q]_{\Z}$.
\end{lem}

\begin{proof} 
Part
 (a) follows from Theorem~\ref{theo2} and Lemma~\ref{lem3}.

To obtain {part} (b), the sequence $0\le\nu_1<\dots<\nu_{n+1}\le q$ is an alternance set of length $n+1$, for example, for the function
\[
f(\nu)=
\begin{cases}
(-1)^{i-1},&\nu=\nu_{i},\ i=1,\dots,n+1,\\
0&\text{otherwise}.
\end{cases}
\]
Note that, by virtue of Theorem~\ref{theo2}, the best uniform approximant is $p^*(\nu)\equiv 0$, and $E(f,L_n)_{\Zinfty}=1$. 
The proof is now complete. 
\end{proof}

\begin{lem}\label{lem6}
Let $p\in L_n$ and $N(p)=n$. There exists a collection of points $\Omega=\{0\le\nu_1<\dots<\nu_{n+1}\le q\}$ such that if for some positive sequence   $\{\rho_i\}_{i=1}^{n+1}$
there holds
\begin{equation}\label{eq2.4}
\sum_{i=1}^{n+1}(-1)^{i}\rho_ip(\nu_i)=0,
\end{equation}
then $p\equiv 0$.
\end{lem}

\begin{proof}
Without loss of generality, we can assume that all $\rho_i=1$. If the polynomial $p$ has
$n$ zeros of the first type, then by the definition of  
$T_0$-systems, $p\equiv 0$.

Let $p$ have a sign change. An interval $[k,k+l]_{\Z}\subset [0,q]_{\Z}$, $l\ge 1$, is called \textit{a closed segment 
 of sign changes} if $p(\nu-1) p(\nu)<0$ for $\nu=k+1,\dots,k+l$ and $p(k-1)p(k)\ge 0$, $p(k+l)p( k+l+1)\ge 0$. It consists of $l+1$ points and contains $l$ zeros of the second type.


Let the polynomial $p$ have $t$ zeros of the first type and $s$  closed segments of sign changes $[k_j,k_j+l_j]_{\Z}$, $j=1,\dots,s$, where $k_{j+1}-k_j-l_j\ge 1$. Then $l_1+\dots+l_s=n-t$.
We form the set $\Omega$ as follows: $\Omega$ consists of
 all $t$ zeros of the first type and $n+1-t$ points from the closed segments 
 $\Omega_j$, $j=1,\dots,s$,
 of sign changes. 
In more detail, set $\Omega_1=[k_1,k_1+l_1]_{\Z}$ 
and construct $\Omega_j$, $j=2,\dots,s$, such that either $\Omega_j=[k_j,k_j+l_j-1]_{\Z}$ or $\Omega_j=[k_j+1 ,k_j+l_j]_{\Z}$. Note that the  number of   points in closed segments is equal to $(l_1+1)+l_2+\dots+l_s=n+1-t$.

If $r_1+i$, $i=0,\dots,l_1$, are the indices of the points $k_1+i$ of the set $\Omega$ in $\Omega_1$
and $p(k_1+i)=(-1)^{i_1 +k_1+i}|p(k_1+i)|$, then $k_1+i=\nu_{r_1+i}$ and
\[
\sum_{i=0}^{l_1}(-1)^{r_1+i}p(\nu_{r_1+i})=(-1)^{i_1+k_1+r_1}\sum_{i=0}^{l_1}|p(\nu_{r_1+i})|.
\]

Let the segments $\Omega_m$, for $m=1,\dots,j-1$, have already been constructed. Let $r_j+i$, $i=0,\dots,l_j-1$, be the indices of points  of the set $\Omega$ in $\Omega_j$ and $p(k_j+i)=(-1)^{i_j+k_j+i}|p(k_j+i)|$, $i=0,\dots,l_j$.

If $(-1)^{i_j+k_j+r_j}=(-1)^{i_1+k_1+r_1}$, then we set $\nu_{r_j+i}=k_j+i$, $i=0, \dots,l_j-1$. We have
\[
\sum_{i=0}^{l_j-1}(-1)^{r_j+i}p(\nu_{r_j+i})=\sum_{i=0}^{l_j-1}(-1)^{r_j+i}p(k_j+i)=(-1)^{i_1+k_1+r_1}\sum_{i=0}^{l_j-1}|p(\nu_{r_j+i})|.
\]

If $(-1)^{i_j+k_j+r_j+1}=(-1)^{i_1+k_1+r_1}$, then we put $\nu_{r_j+i}=k_j+i+1$, $ i=0,\dots,l_j-1$. Again, we have
\[
\sum_{i=0}^{l_j-1}(-1)^{r_j+i}p(\nu_{r_j+i})=\sum_{i=0}^{l_j-1}(-1)^{r_j+i}p(k_j+i+1)=(-1)^{i_1+k_1+r_1}\sum_{i=0}^{l_j-1}|p(\nu_{r_j+i})|.
\]
Thus, the set $\Omega=\{0\le\nu_1<\dots<\nu_{n+1}\le q\}$ is constructed. For this set, in view of \eqref{eq2.4}, there holds
\begin{align*}
0&=\sum_{i=1}^{n+1}(-1)^{i}p(\nu_i)=\sum_{i=0}^{l_1}(-1)^{r_1+i}p(\nu_{r_1+i})+\sum_{j=2}^s\sum_{i=0}^{l_j-1}(-1)^{r_j+i}p(\nu_{r_j+i})\\&
=(-1)^{i_1+k_1+r_1}\Bigl\{\sum_{i=0}^{l_1}|p(\nu_{r_1+i})|+\sum_{j=2}^s\sum_{i=0}^{l_j-1}|p(\nu_{r_j+i})|\Bigr\}.
\end{align*}
Hence, $p(\nu_i)=0$ for $i=1,\dots,n+1$. Then we deduce that $p\equiv 0$, completing the proof. 
\end{proof}

\begin{proof}[Proof of Theorem~\ref{theo3}]
By Lemma~\ref{lem4}, (a) implies (b) and by 
Lemma~\ref{lem5},
items (b) and (c) are equivalent. 
Finally, applying Lemmas~\ref{lem3} and \ref{lem6},  (a) follows from (b). The proof is now complete. 
\end{proof}



\begin{proof}[Proof of Corollary~\ref{cor2}]
To prove this result, it is sufficient to show the following analogue of Lemma~\ref{lem6}: 

{\it If $S^{+}(p)=n$, then there exists a set $\Omega=\{0\le\nu_1<\dots<\nu_{n+1}\le q\}$ such that if for some  positive numbers $\{\rho_i\}_{i=1}^{n+1}$ property \eqref{eq2.4} holds, then $p\equiv 0$.}

The construction of the set $\Omega$ is similar to the one in Lemma~\ref{lem6}, so we will only sketch the proof. As above, we assume that each $\rho_i=1$.

Let us consider four different types of intervals. An interval of the \textit{first} type is a closed segment of sign changes (see the proof of Lemma~\ref{lem6}). An interval of the \textit{second} type is $I=[k,k+l]_{\Z}$, for which
\[
p(k)p(k+l)>0,\quad p(k+1)=\dots=p(k+l-1)=0\quad \text{and}\quad \text{$l-1$ is odd}.
\]
An interval of the \textit{third} type is $I=[k,k+l]_{\Z}$, for which
\[
p(k)p(k+l)<0,\quad p(k+1)=\dots=p(k+l-1)=0\quad \text{and}\quad \text{$l-1$ is even}.
\]

It is worth mentioning that
 intervals of the second and third types contribute to an increase in the number of zeros of the polynomial when zero values are changed to nonzero values, as 
$N(p,[k,k+l]_{\Z})=l-1$ and $S^{+}(p,[k,k+l]_{\Z})=l$.

Intervals of the \textit{fourth} type are $I=[k,k+l]_{\Z}$, consisting of zeros of the first type, for which one of the following 
condition holds:
(i) $k=0$, (ii)  $k+l=q$, (iii) $l+1$ is odd and $p(k-1)p(k+l+1)<0$, (iv) $l+1$ is even and $p(k-1)p(k +l+1)>0$.
Note that changing zero values of $p$ to any nonzero numbers at points of  intervals of the fourth type
cannot increase the number of its zeros.

Denote  the sets of intervals of the first, second, third, and fourth types  by $\Sigma_1$, $\Sigma_2$, $\Sigma_3$, and $\Sigma_4$, respectively. Set $\Sigma=\cup_{i=1}^{3}\Sigma_i$. Note that   intervals from different sets 
$\Sigma_i$
can intersect only at their endpoints. 

We have
\[
n=\sum_{I\in\Sigma_1}(|I|-1)+\sum_{I\in\Sigma_2}(|I|-1)+\sum_{I\in\Sigma_3}(|I|-1)+\sum_{I\in\Sigma_4}|I|.
\]

By $\Omega$ we define  a
set of all\footnote{Here the left or right endpoint of each interval from $\Sigma$  except the first one is excluded, but  each interval from $\Sigma_4$ is included completely; cf.
the proof of Lemma~\ref{lem6}.} points of the intervals from $\Sigma$ 
and $\Sigma_4$.
Then the cardinality of~$\Omega$~is
\[
1+\sum_{I\in\Sigma_1}(|I|-1)+\sum_{I\in\Sigma_2}(|I|-1)+\sum_{I\in\Sigma_3}(|I|-1)+\sum_{I\in\Sigma_4}|I|=n+1.
\]

Let now
$[\nu_{r_1},\nu_{r_1}+l_{r_1}]_{\Z}\subset \Omega$ be the first interval from $\Sigma$ and, for some integer $i_0$,
$(-1)^{i_0+r_1}p(\nu_{r_1})>0$. 
Note that, for $r_1\ge 2$, $0\le\nu_1<\dots<\nu_{r_1-1}$ 
are the points of the set $\Omega\cap \Sigma_4.$
Then $(-1)^{i_0+r_1+i}p(\nu_{r_1}+i)\ge 0$, $i=0,\dots,l_{r_1}$.


Proceeding similarly, let $0\le\nu_1<\dots<\nu_{r_2-1}$ be the points of the set $\Omega$ that have been constructed, and  the second interval $I=[k,k+l]_{\Z}\in \Sigma$ is such that 
$k=\nu_{r_2}$ or $k=\nu_{r_2}-1$.
If $(-1)^{i_0+r_2}p(k)>0$, then we include in $\Omega$ the points of the interval $[\nu_{r_2},\nu_{r_2}+l_{r_2}]_{\Z}$, where $\nu_{r_2}=k$ and $l_{r_2}=l-1$. Otherwise, 
we include in $\Omega$ the points of the interval $[\nu_{r_2},\nu_{r_2}+l_{r_2}]_{\Z}$, where $\nu_{r_2}=k+1$ and $l_{r_2}=l$. In both cases, $|[\nu_{r_2},\nu_{r_2}+l_{r_2}]_{\Z}|=|I|-1$ and $(-1)^{i_0+r_2+i}p(\nu_{r_2}+i)\ge 0$,
$i=0,\dots,l_{r_2}$.

Continuing in the same way, we  complete the construction of the set
$\Omega$.  Then, according to \eqref{eq2.4},
\[
0=\sum_{i=1}^{n+1}(-1)^{i}p(\nu_i)=(-1)^{i_0}\sum_{i=1}^{n+1}|p(\nu_i)|.
\]
From this, $p(\nu_i)=0$ for $i=1,\dots,n+1$. Therefore, $p\equiv 0$, completing the proof.
\end{proof}

The proof of Corollary~\ref{cor2} immediately implies the following result.

\begin{cor}
If the set $\{\varphi_k\}_{k=1}^n\subset C[0,q]_{\Z}$ is a 
$T_{\Z}$-system, $p\in L_n$, and $N(p)=n-1$, then the polynomial $p$ has no intervals of the second and third types.
\end{cor}

To conclude this section, we give an example of a $T_0$-system which is not a  $T_{\Z}$-system.

\begin{exa}\label{example2.7}
For $k=1,\dots,q-1$, define
\[
\varphi_k(\nu)=
\begin{cases}
0,&\nu=0,\dots,k-2,k,\dots,q-1,\\
1,&\nu=k-1,q,
\end{cases}
\quad\varphi_q(\nu)=
\begin{cases}
0,&\nu=0,\dots,q-2,\\
1,&\nu=q-1,\\
-1,&\nu=q.
\end{cases}
\]
Then $\dim{L_q}=q$ and $N_0(\varphi_k)=q-1$, $k=1,\dots,q$. Note that
 the function $\varphi_q(\nu)$ 
 has  a zero of the second type at the point $q$ and hence $N(\varphi_q)=q$.

For an arbitrary polynomial $p(\nu)=\sum_{k=1}^qa_k\varphi_k(\nu)$, the following equalities hold:
\begin{equation}\label{p-nu-a}
p(\nu)=
\begin{cases}
a_{\nu+1},&\nu=0,\dots,q-1,\\
a_1+\dots+a_{q-1}-a_q,&\nu=q.
\end{cases}
\end{equation}
If $p$ has $q$ zeros of the first type, then either all coefficients of the polynomial are zero, or $q-1$ coefficients are zero and $a_1+\dots+a_{q-1}-a_q=0$, that is, in any case all coefficients of $p$ are zeros. Therefore, 
if $p$ is nontrivial, 
$N_0(p)\le q-1$.
Thus, the set $\{\varphi_k\}_{k=1}^q$ is a  $T_0$-system but  not a  $T_{\Z}$-system.

Further, 
define for any  $f\in C[0,q]_{\Z}$ 
\[
p^{*}(\nu)=\sum_{k=1}^{q-1}(f(k-1)+\delta)\varphi_k(\nu)+(f(q-1)-\delta)\varphi_q(\nu),
\]
where
\[
\delta=-\frac{1}{q+1}\Bigl(\,\sum_{k=1}^{q-1}f(k-1)-f(q-1)-f(q)\Bigr).
\]
By \eqref{p-nu-a}, we have
\[
p^{*}(\nu)=
\begin{cases}
f(\nu)+\delta,&\nu=0,\dots,q-2,\\
f(\nu)-\delta,&\nu=q-1,q.
\end{cases}
\]
If $\delta=0$, then $f\equiv p^{*}\in L_q$ and $E_q(f,L_q)_{\infty}=0$. If $\delta\neq 0$, $e(\nu)=p^{*}(\nu)-f(\nu)$,
\[
\varepsilon_{\nu}=
\begin{cases}
\sign \delta,&\nu=0,\dots,q-2,\\
-\sign \delta,&\nu=q-1,q,
\end{cases}
\]
then, by Theorem~\ref{theo2}, for $\nu\in [0,q]_{\Z}$ and all polynomials $p\in L_q$
\[
\varepsilon_{\nu} e(\nu)=\|e\|_{\Zinfty}=|\delta|,\qquad \sum_{\nu=0}^{q}\varepsilon_{\nu} p(\nu)=0.
\]
Consequently, 
$p^{*}\in L_q$ is the unique best uniform approximant of the function $f$ and $E_q(f,L_q)_{\Zinfty}=|\delta|$. Note that the metric projection operator $Pf=p^{*}$ onto the $q$-dimensional subspace $L_q$ is linear (see, e.g., \cite{Bo98}), and {for any function $f\not\in L_q$,} the polynomial $Pf$ does not admit a Chebyshev alternance set of length $q+1$.


\end{exa}

\section{Proofs of Theorems~\ref{theo4}, \ref{theo5}, and \ref{theo6}}\label{sec3}



\begin{proof}[Proof of Theorem~\ref{theo4}] 
Let $k_{l}>0$ be the leading coefficient of the polynomial $P_{l}(\lambda)$. Then, by \eqref{eq1.10}, we have $\beta_{l}k_{l+1}=\rho_{l}k_{l}$. Using this, for the monic polynomials $Q_{l}(\lambda)=k_{l}^{-1}P_{l}(\lambda)$, we obtain the recurrence relation
\[
B_{l}Q_{l-1}(\lambda)+A_{l}Q_{l}(\lambda)+Q_{l+1}(\lambda)=\lambda
Q_{l}(\lambda),
\]
where $B_{0}>0$ is arbitrary and
\begin{equation}\label{A-B}
A_{l}=\frac{\alpha_{l}}{\rho_{l}},\quad
B_{l}=\frac{\gamma_{l-1}\beta_{l-1}}{\rho_{l-1}\rho_{l}}=
\frac{\gamma_{l-1}\rho_{l-1}k_{l-1}^{2}}{\beta_{l-1}\rho_{l}k_{l}^{2}}.
\end{equation}
We have $A_{l}\in \mathbb{R}$ and $B_{l}>0$
for any $l\ge 0$.
These conditions are necessary and sufficient in Favard's theorem \cite[Ch.~II, Th.~6.4]{Ch78} on the existence of a non-decreasing function of bounded variation  with infinite spectrum 
which generates
a positive measure $\mu$ such that
\[
\int_{\mathbb{R}}Q_{l}(\lambda)Q_{m}(\lambda)\,d\mu(\lambda)=B_{0}B_{1}\cdots
B_{l}\delta_{lm}.
\]
Now, to show \eqref{eq1.12}, it is enough to set $B_{0}=\rho_{0}^{-1}$.

According to \cite[Ch.~IV, Th.~2.2]{Ch78}, the measure $\mu$ has a finite support $[a,b]$ if and only if the sequences $A_{l}$ and $B_{l}$ are uniformly bounded. In this case, the endpoints of the interval $[a,b]$ will be the limit points of the zeros of 
$Q_{l}(\lambda)$, and, therefore, of $P_{l}(\lambda)$ as well. By \eqref{A-B}, the boundedness of $A_{l}$ and $B_{l}$ is equivalent to the boundedness of the sequences \eqref{a-g}, which completes the proof. 
\end{proof}

\begin{proof}[Proof of Theorem~\ref{theo5}]
{Recall that 
$\lambda_{l,l}<\dots<\lambda_{1,l}$ 
are the zeros of the polynomials $P_l(\lambda)$, $l\in\mathbb{N}$, 
and 
$\lambda_{q+1}<\dots<\lambda_{1}$ 
are zeros of the polynomial $\widetilde{P}_{q+1}(\lambda)$, cf. \eqref{eq1.13}.}

In what follows, we will use the fact that the zeros of  orthogonal polynomials $P_{l-1}(\lambda)$ and $P_l(\lambda)$, $l\in [2,q]_{\Z}$, as well as $P_q(\lambda)$ and $\widetilde{P}_{q+1}(\lambda)$, interlace (see \cite[Ch.~III, \S\,3.3]{Sz59}):
\begin{equation}\label{eq3.1}
\begin{gathered}
\lambda_{l,l}<\lambda_{l-1,l-1}<\lambda_{l-1,l}<\lambda_{l-2,l-1}<\dots<
\lambda_{2,l-1}<\lambda_{2,l}<\lambda_{1,l-1}<\lambda_{1,l},\\
\lambda_{q+1}<\lambda_{q,q}<\lambda_{q}<\lambda_{q-1,q}<\dots<
\lambda_{2,q}<\lambda_{2}<\lambda_{1,q}<\lambda_{1}.
\end{gathered}
\end{equation}

Let the interval $(a,b)$ contain the zeros of the polynomial $P_q(\lambda)$  and let $k\in [1,q+1]_{\Z}$.  If $s=0$, then $\psi_{k}(0)=P_0(\lambda_{k})=1$, so further $s\in [1,q]_{\Z}$.

We will also use the following observation. If $s\in [1,q]_{\Z}$, $l\in [1,s]_{\Z}$, $\lambda_{s+1,s}=a$, $\lambda_ {0,s}=b$, and $\lambda_{l,s}<\lambda_{k}<\lambda_{l-1,s}$, then
\begin{equation}\label{eq3.2}
\sign \psi_k(s)=\sign P_s(\lambda_{k})=(-1)^{l-1},
\end{equation}
and if $\lambda_{k}=\lambda_{l,s}$ or $\lambda_{k}=\lambda_{l-1,s}$, then
\begin{equation}\label{eq3.3}
\psi_k(s)=P_s(\lambda_{k})=P_s(\lambda_{l,s})=P_s(\lambda_{l-1,s})=0.
\end{equation}

Further, for $k=1$ and $s\in [1,q]_{\Z}$, the inequality $\lambda_{1,s}<\lambda_{1}$ yields  
\begin{equation}\label{psi-1-pos}
\psi_1(s)>0,\quad s\in [0,q]_{\Z}. 
\end{equation}
Therefore, $N(\psi_1)=0$. For $k=q+1$ and $s\in [1,q]_{\Z}$, from the inequalities \eqref{eq3.1} we have $\lambda_{q+1}<\lambda_{s,s}$, so $\sign \psi_{q+1}(s)=(-1)^{s}$ and $\psi_{q+1}(0)=1$. Thus, $N(\psi_{q+1})=q$.

Let $k\in [2,q]_{\Z}$. Since, according to \eqref{eq3.1}, $\lambda_{k}<\lambda_{k-1,q}$, there exists a smallest $i_1\le q$ for which
\begin{equation}\label{eq3.4}
\lambda_{k,i_1}<\lambda_{k-1,i_1-1}\le\lambda_{k}<\lambda_{k-1,i_1},\quad i_1\ge k.
\end{equation}
If $i_1\le s\le q$, then, by \eqref{eq3.1} and \eqref{eq3.4},
\[
\lambda_{k,i_1}\le\lambda_{k,s}\le \lambda_{k,q}< \lambda_{k}<\lambda_{k-1,i_1}\le\lambda_{k-1,s}\le\lambda_{k-1,q}.
\]
Using \eqref{eq3.2}, we then derive
\begin{equation}\label{eq3.5}
\sign \psi_k(s)=(-1)^{k-1}.
\end{equation}

Since, according to \eqref{eq3.1}, $\lambda_{k}<\lambda_{k-1,i_1}<\lambda_{k-2,i_1-1}$, there exists a smallest $i_2\le i_1-1 $ for which
\begin{equation}\label{eq3.6}
\lambda_{k-1,i_2}<\lambda_{k-2,i_2-1}\le\lambda_{k}<\lambda_{k-2,i_2}, \quad i_2\ge k-1.
\end{equation}
If $i_2\le s<i_1-1$, then, by \eqref{eq3.1}, \eqref{eq3.4} and \eqref{eq3.6},
\[
\lambda_{k-1,i_2}\le\lambda_{k-1,s}<\lambda_{k-1,i_1-1}\le\lambda_{k}<\lambda_{k-2,i_2}\le\lambda_{k-2,s}<\lambda_{k-2,i_1-1}.
\]
Therefore, using \eqref{eq3.2} and \eqref{eq3.3}, we obtain
\[
\sign \psi_k(s)=(-1)^{k-2},\quad \sign \psi_k(i_1-1)=(-1)^{k-2}\vee 0.
\]

Proceeding similarly for $j\in [2,k-1]_{\Z}$, we construct points $\{i_j\}_{j=1}^{k-1}$ such that
\[
k-j+1\le i_j\le i_{j-1}-1, \quad \lambda_{k-j+1,i_j}<\lambda_{k-j,i_j-1}\le\lambda_{k}<\lambda_{k-j,i_j},
\]
and for $i_j\le s<i_{j-1}-1$
\[
\lambda_{k+1-j,i_j}\le\lambda_{k+1-j,s}<\lambda_{k+1-j,i_{j-1}-1}\le\lambda_{k}<\lambda_{k-j,i_j}\le\lambda_{k-j,s}<\lambda_{k-j,i_{j-1}-1},
\]
\begin{equation}\label{eq3.7}
\sign \psi_k(s)=(-1)^{k-j},\quad \sign \psi_k(i_{j-1}-1)=(-1)^{k-j}\vee 0.
\end{equation}
In the case $1\le s<i_{k-1}-1$, we have $\lambda_{1,s}<\lambda_{1,i_{k-1}-1}\le\lambda_{k}$. Hence, in light of \eqref{eq3.1},
\begin{equation}\label{eq3.8}
\sign \psi_k(s)=1,\quad \sign \psi_k(i_{k-1}-1)=1\vee 0.
\end{equation}
If for all $j\in[2, k-1]_{\Z}$, $\lambda_{k-j,i_j-1}\neq\lambda_{k}$, then, according to \eqref{eq3.5}, \eqref{eq3.7}, \eqref{eq3.8}, only the points $\nu_j=i_j$, $j\in [1,k-1]_{\Z}$, are zeros of $\psi_k(\nu)$ and, moreover, $N(\psi_k)=k-1$.

Let $\psi_k(i_j-1)=0$. Since $\psi_k(i_{j+1})\neq 0$, it follows that $i_{j+1}\le i_j-2<i_j-1$ and
\[
\sign{}(\psi_k(i_{j}-2)\psi_k(i_{j}))=(-1)^{k-j-1}(-1)^{k-j}=-1.
\]
Therefore, the zero $\nu_j=i_j$ is replaced by $\nu_j=i_j-1$ and we again have $N(\psi_k)=k-1$. Along the way, we proved that if  $\psi_{k}(\nu)=0$,
$k\in [2,q]_{\Z}$,
then $\nu\in [1,q-1]_ {\Z}$ and $\psi_{k}(\nu-1)\psi_{k}(\nu+1)<0$. This implies that $S^{-}(\psi_k)=S^{+}(\psi_k)=k-1$, completing the proof.
\end{proof}

\subsection{Proof of Theorem~\ref{theo6}}

\textit{Step 1.} First, let us 
 rewrite the discrete Sturm-Liouville problems \eqref{eq1.10} and  \eqref{eq1.15} and then derive the Christoffel-Darboux formulas. 

The self-adjoint form of 
 the problem \eqref{eq1.10} is given by
\begin{equation*}
\nabla(w_l\Delta P_l(\lambda))+c_lP_l(\lambda)=\lambda d_lP_l(\lambda), \quad l\in \mathbb{Z}_{+},\quad \lambda\in\mathbb{C},
\end{equation*}
where $\Delta b_{l}=b_{l+1}-b_{l}$, $\nabla b_{l}=b_{l}-b_{l-1}$ are the forward and backward differences, respectively, $d_l$ is defined in \eqref{eq1.12}, and
\begin{equation}\label{eq3.10}
w_l=\frac{\beta_ld_l}{\rho_l},\quad c_l=\frac{(\alpha_l+\beta_l+\gamma_{l-1})d_l}{\rho_l}.
\end{equation}
Indeed, using \eqref{eq1.10}, \eqref{eq1.12} and \eqref{eq3.10}, we obtain
\[
w_{l-1}=\frac{\gamma_{l-1}d_l}{\rho_l},\quad \nabla(w_l\Delta P_l(\lambda))=w_lP_{l+1}(\lambda)+w_{l-1}P_{l-1}(\lambda)-(w_{l}+w_{l-1})P_{l}(\lambda),
\]
and
\[
\nabla(w_l\Delta P_l(\lambda))+c_lP_l(\lambda)-\lambda d_lP_l(\lambda)
=\frac{d_l}{\rho_l}\,\bigl\{\beta_lP_{l+1}(\lambda)+\alpha_lP_{l}(\lambda)+\gamma_{l-1}P_{l-1}(\lambda)-\lambda\rho_lP_{l}(\lambda)\bigr\}.
\]
Now we can easily derive  the Christoffel-Darboux formula
\begin{equation}\label{eq3.11}
\sum_{s=0}^{l}d_{s}P_{s}(x)P_{s}(y)=w_{l}\,\frac{P_{l+1}(x)P_{l}(y)-P_{l}(x)P_{l+1}(y)}{x-y},\quad x\neq y.
\end{equation}
To see this, using the Lagrange formula  \cite{ABGR05}
\[
(x-y)d_sP_s(x)P_s(y)=P_s(y)\nabla(w_s\Delta P_s(x))-P_s(x)\nabla(w_s\Delta
P_s(y)),
\]
we have 
\begin{align*}&
(x-y)\sum_{s=0}^{l}d_sP_s(x)P_s(y)=\sum_{s=0}^{l}\bigl\{P_s(y)\nabla(w_s\Delta P_s(x))-P_s(x)\nabla(w_s\Delta P_s(y))\bigr\}\\
&=\sum_{s=0}^{l}\bigl\{w_sP_{s+1}(x)P_s(y)+w_{s-1}P_s(y)P_{s-1}(x)-w_sP_{s+1}(y)P_s(x)-w_{s-1}P_{s}(x)P_{s-1}(y)\bigr\}\\
&=w_l\bigl\{P_{l+1}(x)P_l(y)-P_{l}(x)P_{l+1}(y)\bigr\}.
\end{align*}
Further, by \eqref{eq3.11}, it follows that $w_{l}=k_{l}d_{l}/k_{l+1}$, where $k_{l}$ is the leading coefficient of the polynomial $P_{l }$. Applying \eqref{eq1.12} and \eqref{eq3.10}, we get
\[
k_{0}=1,\quad k_{l}=\frac{\rho_0\cdots\rho_{l-1}}{\beta_0\cdots\beta_{l-1}}, \quad l\ge 1.
\]

Moreover, the Sturm-Liouville problem \eqref{eq1.15} can also be written in an equivalent self-adjoint form as follows:
\begin{equation}\label{eq3.12}
\begin{gathered}
\nabla(w_l\Delta \psi(\nu))+c_l\psi(\nu)=\lambda d_l\psi(\nu),\\
\psi(-1)=0,\quad\psi(q+1)-\eta\psi(q)=0.
\end{gathered}
\end{equation}
From \eqref{eq3.12}, as in the proof of \eqref{eq3.11}, we obtain  the following Lagrange equality:
\begin{equation*}
d_{\nu}(\lambda_l-\lambda_k)\psi_l(\nu)\psi_k(\nu)=
\nabla(w_{\nu}\{\psi_l(\nu)\Delta \psi_k(\nu)-\psi_k(\nu)\Delta \psi_l(\nu)\}).
\end{equation*}
It is important to mention that  the Christoffel-Darboux formula 
for the eigenfunctions  $\psi_k$
has the form
\begin{equation}\label{psi-kd}
\sum_{\nu=0}^{l}d_{\nu}\psi_k(\nu)\psi_s(\nu)=w_{l}\,\frac{\psi_k(l+1)\psi_s(l)-\psi_k(l)\psi_s(l+1)}{\lambda_k-\lambda_s},
\end{equation}
which immediately implies the discrete orthogonality of the eigenfunctions. Indeed, for $k\neq s$, due to the boundary condition in \eqref{eq1.15},
\begin{align*}
(\psi_k,\psi_s)&=\sum_{\nu=0}^{q}d_{\nu}\psi_k(\nu)\psi_s(\nu)=w_{q}\,\frac{\psi_k(q+1)\psi_s(q)-\psi_k(q)\psi_s(q+1)}{\lambda_k-\lambda_s}\\
&=w_{q}\,\frac{(\psi_k(q+1)-\eta\psi_k(q))\psi_s(q)-(\psi_s(q+1)-\eta\psi_s(q))\psi_k(q)}{\lambda_k-\lambda_s}=0.
\end{align*}


\textit{Step~2.} Now we will prove several lemmas on zeros and sign changes of discrete functions.

\begin{lem}\label{lem7}
If $f\in C[0,q]_{\Z}$, then
\begin{equation}\label{eq3.15}
S^{-}(f)\le N(f)\le S^{+}(f).
\end{equation}
\end{lem}

\begin{proof} 
It suffices to compare $S^{-}(f,[m,n]_{\Z})$, $N(f,[m,n]_{\Z})$, and $S ^{+}(f,[m,n]_{\Z})$  at interior points of 
the intervals $[m,n]_{\Z}$,
where $f$ equals zero. We may assume that $n-m\ge 2$.



\begin{enumerate}
\item[\bigcdot] If $m=0$, $f(0)=f(1)=\dots=f(n-1)=0$, and $f(n)\neq 0$, then $S^{-}(f ,[0,n]_{\Z})=0$, $N(f,[0,n]_{\Z})=n$, $S^{+}(f,[0,n] _{\Z})=n$ and  thus \eqref{eq3.15} is valid.

\item[\bigcdot] If $n=q$, $f(m)\neq 0$, $f(m+1)=\dots=f(q)=0$, and $f(n)\neq 0$, then $S^ {-}(f,[m,q]_{\Z})=0$, $N(f,[m,q]_{\Z})=q-m$, $S^{+}(f, [m,q]_{\Z})=q-m$ and \eqref{eq3.15} follows.

\item[\bigcdot] If $f(m)f(n)<0$ (for example, $f(m)>0$, $f(m+1)=\dots=f(n-1)=0$, $f(n )<0$), then $S^{-}(f,[m,q]_{\Z})=1$, $N(f,[m,q]_{\Z})=n-m-1 $, $S^{+}(f,[m,n]_{\Z})=n-m-1$ for even $n-m$ and $S^{+}(f,[m,n]_{\Z})=n-m$ for odd $n-m$. Hence, \eqref{eq3.15} is valid.
  
\item[\bigcdot] If $f(m)f(n)>0$ (for example, $f(m)>0$, $f(m+1)=\dots=f(n-1)=0$, $f(n )>0$), then $S^{-}(f,[m,q]_{\Z})=1$, $N(f,[m,q]_{\Z})=n-m-1 $, $S^{+}(f,[m,n]_{\Z})=n-m-1$ for odd $n-m$ and $S^{+}(f,[m,n]_{\Z})=n-m$ for even $n-m$, implying  \eqref{eq3.15}.
\end{enumerate}



The proof is now complete.
\end{proof}

The next few results correspond to discrete analogues of Rolle's theorem.
Below, for integers $m<n$, we say $(m,n)$ is a pair of  \textit{generalized 
sign change} of $f$ if $f(m)f(n)<0$ and
$f(m+1)=\ldots=f(n-1)=0$.
Clearly, an ordinary sign-change pair is also generalized  (for $n-m=1$).
Note that the number of pairs of generalized sign changes of $f$ on
$[a,b]_{\Z}$ is exactly $S^{-}(f,[a,b]_{\Z})$.

\begin{lem}\label{lem8}
If $m<n$ and $f(m)f(n)<0$, then the function $f$ has a zero on $[m+1,n]_{\Z}$ and a generalized sign change  on $[m,n]_{\Z}$.
\end{lem}

\begin{proof} 
We can assume that $f(m)<0$ and $f(n)>0$. If $n=m+1$, then $n$ is a zero of the second type and  $(m,n)$ is the pair of a sign change. Let $n>m+1$ and $s_1\le n$ be the largest integer such that $f(m)<0,f(m+1)<0,\dots,f(s_1-1)<0 $. Then either $f(s_1)>0$ or $f(s_1)=0$. If $f(s_1)>0$, then $s_1\ge m+1$ is a zero of the second type and  $(s_1-1,s_1)$ is the pair of a sign change.
Suppose $f(s_1)=0$ and $s_2<n$ is the largest integer such that $f(s_1)=f(s_1+1)=\dots=f(s_2-1)=0$. Then either $f(s_2)>0$ or $f(s_2)<0$. If $f(s_2)>0$, then $s_1$, $s_1+1$,\dots, $s_2-1$ are zeros of the first type and $(s_1-1,s_2)$ 
is the pair of a generalized sign change.
 If $f(s_2)<0$, using similar arguments, we arrive at the required statement. 
\end{proof}

\begin{lem}\label{lem9}
Let $m<n$ be zeros of $f$. Then
the differences $\nabla f$ on $[m+1,n]_{\Z}$ and $\Delta f$ on $[m,n-1 ]_{\Z}$ have zeros.
\end{lem}

\begin{proof} 
First, let us prove the lemma for $\nabla f$. Consider the possible values  $f$ at the zeros $m$ and $n$. 
We will repeatedly use the trivial formula
\begin{equation}\label{eq3.16}
f(l)-f(k)=\sum_{\nu=k+1}^{l}\nabla f(\nu).
\end{equation}

\begin{enumerate}
\item[\bigcdot] If $f(m)=f(n)=0$, then according to \eqref{eq3.16}, either $\nabla f (\nu)=0$ for all $\nu\in[m+1,n]_{\Z}$, or $\nabla f(\nu_1)\nabla f(\nu_2)<0 $ for some $\nu_1,\nu_2\in[m+1,n]_{\Z}$, and by Lemma~\ref{lem8} $\nabla f$  has a zero on $[m+1,n]_{\Z}$.

\item[\bigcdot] If $f(m)=0$, $f(n-1)>0$, $f(n)<0$, then $f(n-1)-f(m)>0$, $\nabla f(n)<0$. According to \eqref{eq3.16}, for some $s\in[m+1,n-1]_{\Z}$, $\nabla f(s)>0$, so $\nabla f(s)\nabla f(n)<0$. By Lemma~\ref{lem8}, $\nabla f$ on $[m+1,n]_{\Z}$ has a zero. The case $f(m)=0$, $f(n-1)<0$, $f(n)>0$ is considered similarly.

\item[\bigcdot] If $f(m-1)>0$, $f(m)<0$, $f(n)=0$, then $\nabla f(m)<0$, $f(n)-f( m)>0$. According to \eqref{eq3.16}, for some $s\in[m+1,n]_{\Z}$, $\nabla f(s)>0$, therefore $\nabla f(m)\nabla f (s)<0$. In light of Lemma~\ref{lem8}, $\nabla f$ on $[m+1,n]_{\Z}$ has a zero. The case $f(m-1)<0$, $f(m)>0$, $f(n)=0$ is treated similarly.

\item[\bigcdot] If $f(m-1)>0$, $f(m)<0$, $f(n-1)>0$, $f(n)<0$, then $\nabla f(m)< 0$, $f(n-1)-f(m)>0$. According to \eqref{eq3.16}, for some $s\in[m+1,n-1]_{\Z}$, $\nabla f(s)>0$, so $\nabla f(m)\nabla f(s)<0$. By Lemma~\ref{lem8}, $\nabla f$ on $[m+1,n]_{\Z}$ has a zero. 

\item[\bigcdot] If $f(m-1)>0$, $f(m)<0$, $f(n-1)<0$, $f(n)>0$, then $\nabla f(m)< 0$, $\nabla f(n)<0$. Using again Lemma~\ref{lem8}, $\nabla f$ on $[m+1,n]_{\Z}$ has a zero. The remaining two cases are considered similarly.
\end{enumerate}

Finally,  
using the equality
$
f(l)-f(k)=\sum_{\nu=k}^{l-1}\Delta f(\nu)
$
one similarly proves
the statement for $\Delta f$.
The proof is now complete.
\end{proof}

\begin{lem}\label{lem10}
Let {$f\in C[0,q]_{\Z}$, $0\le m<n\le m'<n'\le q$}, and the pairs $(m,n)$ and $(m',n')$ be consecutive generalized sign changes of $f$. Then the differences $\nabla f$ on $[m+1,m'+1]_{\Z}$ and $\Delta f$ on $[n-1,n'-1]_{\Z}$ have generalized  sign  changes. 
\end{lem}

\begin{proof} 
Let $f(m)<0$ and $f(n)>0$. Then
\[
f(\nu)=0,\quad \nu\in [m+1,n-1]_{\Z},\quad f(\nu)>0,\quad \nu\in [n,m']_{\Z},
\]
\[
f(\nu)=0,\quad \nu\in [m'+1,n'-1]_{\Z},\quad f(n')<0.
\]
Since $\nabla f(m+1)=f(m+1)-f(m)>0$ and $\nabla f(m'+1)=f(m'+1)-f(m')<0 $, by Lemma~\ref{lem8},
$\nabla f$ has a generalized sign change on $[m+1,m'+1]_{\Z}$. The case $f(m)>0$ and $f(n)<0$ and 
the fact that the difference $\Delta f$ has a generalized  sign change on $[n-1,n'-1]_{\Z}$  can be obtain similarly.
\end{proof}

\begin{lem}\label{lem10'}
Let $f\in C[0,q]_{\Z}$, 
$f(-1)=f(q+1)=0$,
$0\le m<n\le q$,
and the pair $(m,n)$ be a generalized sign change of $f$. Then the difference $\nabla f$ on $[0,m+1]_{\Z}$ and $[n,q+1]_{\Z}$ has generalized sign changes.
\end{lem}

\begin{proof} 
Let $f(\nu)=0$ for $\nu\in [-1,m-1]_{\Z}$. If $f(m)>0$ and $f(n)<0$, then $\nabla f(m)>0$. Since $\nabla f(m+1)<0$, $(m,m+1)$ 
is the pair of a sign change
 of $\nabla f$. The case $f(m)<0$ and $f(n)>0$ is treated similarly.

Suppose the function $f(\nu)$ has a nonzero value on $[0,m-1]_{\Z}$, $m\ge 1$, for example, $f(0)>0$.
If $f(m)>0$ and $f(n)<0$, then $\nabla f(0)>0$ and $\nabla f(m+1)<0$. By Lemma~\ref{lem8}, $\nabla f$  has a generalized  sign change on $[0,m+1]_{\Z}$. If $f(m)<0$ and $f(n)>0$, then $f(m)-f(0)<0$. Then for some $\nu_1\in [1,m]_{\Z}$, $\nabla f(\nu_1)<0$
and $\nabla f(0)>0$. By Lemma~\ref{lem8}, $\nabla f$  has a generalized sign change on $[0,m]_{\Z}$. The case $f(0)<0$ is treated similarly. 
Therefore, $\nabla f$ on $[0,m+1]_{\Z}$ always has a generalized sign change.

The fact that  the difference $\nabla f$ 
has a generalized sign change  on $[n,q+1]_{\Z}$ can be shown similarly. 
The proof is now complete.
\end{proof}

\begin{defn}
By $K(\Cdot)$ we denote any one of the following quantities: $N(\Cdot)$, $S^{-}(\Cdot)$, or $S^{+}(\Cdot)$.
\end{defn}

We will use the simple facts that 
$K(f,[a,b]_{\Z})=K(f,[a,c]_{\Z})+K(f,[c,b]_{\Z})$ 
whenever $a\le c\le b$ with $f(c)\ne 0$,
and that $K(fg)=K(f)$ if $g$ preserves the sign.

\begin{proposition}\label{lem11}
If $f\in C[0,q]_{\Z}$, $f\ne 0$, then
\begin{equation}\label{eq3.17}
K(\nabla f,[1,q]_{\Z})\ge K(f,[0,q]_{\Z})-1,
\end{equation}
\begin{equation}\label{eq3.18}
K(\Delta f,[0,q-1]_{\Z})\ge K(f,[0,q]_{\Z})-1.
\end{equation}
If, in addition, $f(-1)=f(q+1)=0$, then
\begin{equation}\label{eq3.17'}
K(\nabla f,[0,q+1]_{\Z})\ge K(f,[0,q]_{\Z})+1.
\end{equation}
\end{proposition}

\begin{rem}
Note that 
\eqref{eq3.18} 
was proved  in
\cite[Proposition~5.1]{Ha78} (see also 
the monograph
\cite[Theorem~1.8.1]{Ag00}) for 
$K(\Cdot)=N(\Cdot)$.
\end{rem}

\begin{proof}[Proof of Proposition \ref{lem11}] 
As usual, we denote 
$K(f)=
K(f,[0,q]_{\Z})
.$
{The cases $f>0$
and $f<0$ 
 on $[0,q]_{\Z}$ are straightforward and may be omitted.}

(a) The cases $K(\Cdot)=N(\Cdot)$ or $K(\Cdot)=S^{-}(\Cdot)$.
To  prove inequality \eqref{eq3.17} in the first case,
let $\{\nu_1<\dots<\nu_{s}\}\subset [0,q]_{\Z}$ be the zeros of  $f$. By Lemma~\ref{lem9} on the intervals $[\nu_1+1,\nu_2]_{\Z},\dots,[\nu_{s-1}+1,\nu_{s}]_{\Z}$, the function $\nabla f$ has at least $s-1$ zeros. Therefore, $N(\nabla f,[1,q]_{\Z})\ge N(f)-1$.

In the second case, 
 let $\{(m_i,n_i)\}_{i=1}^{s}\subset [0,q]_{\Z}\times
[0,q]_{\Z}$, $m_i<n_i\le m_{i+1}<n_{i+1}$, be the generalized sign changes of  $f$ on $[0,q]_{\Z}$. By Lemma~\ref{lem10}, $\nabla f$ has generalized sign changes on the intervals $[m_1+1,m_2+1]_{\Z},\dots,[m_{s-1}+1,m_s+1]_{\Z}$. Since this number is not less than $s-1$, $S^{-}(\nabla f,[1,q]_{\Z})\ge S^{-}(f)-1$.

{Inequality} \eqref{eq3.18} can be established similarly.  

Suppose $f(-1)=f(q+1)=0$. By Lemma~\ref{lem9},
the difference $\nabla f$ has at least two additional zeros on $[0,\nu_{1}]_{\Z}$ and $[\nu_{s}+1,q+1]_{\Z}$.
Similarly, by Lemma~\ref{lem10'},
$\nabla f$ has at least two generalized sign changes on the intervals $[0,m_{1}+1]_{\Z}$ and $[\nu_{s},q+1]_{\Z}$.
This establishes \eqref{eq3.17'}.

\medbreak
(b) The case $K(\Cdot)=S^{+}(\Cdot)$.
For the sake of brevity, we will prove 
the intermediate inequality 
\begin{equation}\label{Snabla-S}
S^{+}(\nabla f)\ge S^{+}(f),\quad f(-1)=0,
\end{equation}
by induction on $q$. Inequalities \eqref{eq3.17} and \eqref{eq3.17'} can be proved in a similar way.

Let $n=q=1$ and $f(0)=\nabla f(0)>0$. If $f(1)>0$, then $S^{+}(\nabla f)=0$ and  \eqref{Snabla-S} is valid. If $f(1)\le 0$, then $\nabla f(1)<0$, so $S^{+}(\nabla f)=S^{+}(f)=1$ and \eqref {Snabla-S} holds. The case $f(0)<0$ can be treated similarly. Let $f(0)=\nabla f(0)=0$. If $f(1)\neq 0$, then $S^{+}(\nabla f)=S^{+}(f)=1$ and \eqref{Snabla-S} is valid.


Assume that inequality \eqref{Snabla-S} holds for all $n\le q-1$. To prove it for  $n=q$,
if $n\le q-1$ is chosen so that $f(n)\neq 0$, $\nabla f(n)\neq 0$ and $S^{+}(\nabla f,[n,q ]_{\Z})\ge S^{+}(f,[n,q]_{\Z})$, then by the inductive assumption we obtain \eqref{Snabla-S} as follows:
\[
S^{+}(\nabla f)=S^{+}(\nabla f,[0,n]_{\Z})+S^{+}(\nabla f,[n,q]_{\Z})\ge S^{+}(f,[0,n]_{\Z})+S^{+}(f,[n,q]_{\Z})=S^{+}(f).
\]

\begin{enumerate}
\item[\bigcdot]
Let $f(q-1)>0$. If $f(q)>0$, then, by the inductive assumption,
\[
S^{+}(\nabla f)\ge S^{+}(\nabla f,[0,q-1]_{\Z})\ge S^{+}(f,[0,q-1]_{\Z})=S^{+}(f).
\]

\item[\bigcdot]
Let $f(q-1)>0$, $f(q) \le 0$, and $s \le q+1$ be  such that
\[
f(q-s)\le 0,\ f(q-s+1)>0,\ \dots,\ f(q-1)>0.
\]
\end{enumerate}
Since $\nabla f(q-s+1)>0$ and $\nabla f(q)<0$,  Lemma~\ref{lem8} implies that  $S^{+}(\nabla f,[q-s+1,q]_{\Z}) \ge 1$ and
\[
S^{+}(\nabla f,[q-s+1,q]_{\Z})\ge 1=S^{+}(f,[q-s+1,q]_{\Z}).
\]
The case  $f(q-1)<0$ can be treated similarly.

Assume now that 
$f(q-1)=\dots=f(q-s_1+1)=0$ and
\[
\ f(q-s_1)>0,\ \dots,\ f(q-s_1-s_2+1)>0,\ f(q-s_1-s_2)\le 0.
\]

\begin{enumerate}
\item[\bigcdot] If $f(q)=0$, then
\[
S^{+}(\nabla f,[q-s_1-s_2+1,q]_{\Z})\ge S^{+}(f,[q-s_1-s_2+1,q]_{\Z})=s_1.
\]

\item[\bigcdot] If $f(q)>0$, then
\[
S^{+}(\nabla f,[q-s_1-s_2+1,q]_{\Z})\ge S^{+}(f,[q-s_1-s_2+1,q]_{\Z})=
\begin{cases}
s_1-1,&\text{$s_1$ odd},\\
s_1,&\text{$s_1$ even}.
\end{cases}
\]

\item[\bigcdot] If $f(q)<0$, then
\[
S^{+}(\nabla f,[q-s_1-s_2+1,q]_{\Z})\ge S^{+}(f,[q-s_1-s_2+1,q]_{\Z})=
\begin{cases}
s_1,&\text{$s_1$ odd},\\
s_1-1,&\text{$s_1$ even}.
\end{cases}
\]
\end{enumerate}
The remaining cases can be reduced to the ones above
by replacing $f$ with $-f$. 
The proof is now complete.
%
\end{proof}

In order to prove Theorem~\ref{theo6}, we employ the discrete version of the Liouville method (see, e.g., \cite{BH20} for the continuous case). 
In addition to the polynomial $
V(\nu)=\sum_{k=m}^{n}a_k\psi_k(\nu)$
(see \eqref{eq1.17}), we define 
\begin{equation*}
V_r(\nu)=\sum_{s=m}^{n}(\lambda_1-\lambda_s)^ra_s\psi_s(\nu), \quad r\in\mathbb{Z}.
\end{equation*}

\begin{lem}\label{lem13}
We have
\[
K(V_1)\ge K(V).
\]
\end{lem}

\begin{proof}
Multiplying both sides of \eqref{psi-kd} 
with $k=1$
by $(\lambda_1-\lambda_s)a_s$
  and summing  over $s$, for $l\in [ 0,q-1]_{\Z}$, we derive
\begin{align}
g(l)&=\sum_{\nu=0}^ld_{\nu}\psi_1(\nu)V_{1}(\nu)
=w_{l}\{\psi_1(l+1)V(l)-\psi_1(l)V(l+1)\}\notag\\ \label{eq3.22}
&=-w_{l}\psi_1(l)\psi_1(l+1)\Bigl\{\frac{V(l+1)}{\psi_1(l+1)}-\frac{V(l)}{\psi_1(l)}\Bigr\}=
-w_{l}\psi_1(l)\psi_1(l+1)\Delta\Bigl(\frac{V(l)}{\psi_1(l)}\Bigr).
\end{align}
Due to the boundary condition in \eqref{eq1.15}, we have
$g(-1)=0$ and
\begin{align*}
g(q)&=w_{q}\{\psi_1(q+1)V(q)-\psi_1(q)V(q+1)\}\\&=w_{q}\{(\psi_1(q+1)-\eta \psi_1(q))V(q)-\psi_1(q)(V(q+1)-\eta V(q))\}=0.
\end{align*}

Applying Proposition~\ref{lem11}, formula \eqref{eq3.22}, and taking into account the positivity of $\psi_1(l)$ (see~\eqref{psi-1-pos}), $d_l$, $w_l$, we obtain 
\begin{align*}
K(V_1)&=K\Bigl(\frac{\nabla g(l)}{d_l\psi_1(l)}\Bigr)=K(\nabla g,[0,q]_{\Z})\ge K(g,[0,q-1]_{\Z})+1\\
&=K\Bigl(-w_s\psi_1(l)\psi_1(l+1)\Delta\Bigl(\frac{V(l)}{\psi_1(l)}\Bigr),[0,q-1]_{\Z}\Bigr)+1\\
&=K\Bigl(\Delta\Bigl(\frac{V(l)}{\psi_1(l)}\Bigr),[0,q-1]_{\Z}\Bigr)+1\ge K\Bigl(\frac{V(l)}{\psi_1(l)},[0,q]_{\Z}\Bigr)=K(V),
\end{align*}
completing the proof of the lemma.
\end{proof}


\begin{proof}[\textit{Step~3 of the proof of Theorem~\ref{theo6}.}]
Note that in the case $m=n=1$ Theorem~\ref{theo6} holds since $\psi_1$ is positive. Thus, we may assume that $n\ge 2$.
Moreover,  the lower bound in \eqref{eq1.17+} for $m=1$ holds trivially. 

Applying Lemma~\ref{lem13} repeatedly, we obtain
\[
K(V_r/(\lambda_{1}-\lambda_n)^r)=K(V_r)\ge K(V),\quad r\ge 1,\quad m\ge 1,
\]
\[
K(V_r/(\lambda_{1}-\lambda_m)^r)=K(V_r)\le K(V),\quad r\le-1,\quad m\ge 2.
\]
Since
\[
\lim_{r\to +\infty}\|\psi_{n}-V_r/(\lambda_{1}-\lambda_n)^r)\|_{\Zinfty}=
\lim_{r\to -\infty}\|\psi_{m}-V_r/(\lambda_{1}-\lambda_m)^r)\|_{\Zinfty}=0,
\]
by Theorem~\ref{theo5} for sufficiently large $N>0$ we deduce 
\[
K(\psi_{n})=K(V_r/(\lambda_{1}-\lambda_n)^r),\quad r>N,
\]
\[
K(\psi_{m})=K(V_r/(\lambda_{1}-\lambda_m)^r),\quad r<-N.
\]
Thus,\[
n-1=K(\psi_{n})=\lim_{r\to+\infty}K(V_r/(\lambda_{1}-\lambda_n)^r)=\lim_{r\to+\infty}K(V_r)\ge K(V),
\]
\[
m-1=K(\psi_{m})=\lim_{r\to -\infty}K(V_r/(\lambda_{1}-\lambda_m)^r)=\lim_{r\to-\infty}K(V_r)\le K(V).
\]
This and Lemma~\ref{lem7} complete the proof.
\end{proof}

\section{Proof of Theorem~\ref{theo7}}\label{sec4}

Let $\nu\in [0,q]_{\Z}$,
 $1\le\nu_1<\dots<\nu_{m}\le q$,
 and 
\[
R_{m}(\nu)=R_{m}(\nu,\nu_1,\dots,\nu_{m})=\prod_{j=1}^{m}(\nu_j-\nu).
\]
Let also
 $\{\varphi_k\}_{k=1} ^{m+1}$ be a 
$T_{\Z}$-system on $[0,q]_{\Z}$ and
\begin{equation}\label{eq4.1}
D_{m+1}(\nu)=D_{m+1}(\nu,\nu_1,\dots,\nu_{m})
=\Delta\!\left(\setlength{\arraycolsep}{3pt}
\begin{matrix}
\varphi_1, & \dots, &\varphi_{m},& \varphi_{m+1}\\
\nu_{m}, & \dots, &\nu_1, &\nu
\end{matrix}\right).
\end{equation}
Observe that  $\nu_j$, $j=1,\dots,m$, are the only zeros of the polynomial $D_{m+1}$.
By applying Theorem~\ref{theo3} 
and basic properties of  determinants, 
we obtain the precise  distribution of the signs of $D_{m+1}(\nu)$.

\begin{lem}\label{lem14} 
For all $\nu\in [0,q]_{\Z}$,
\begin{equation*}
\sign D_{m+1}(\nu)=\sign D_{m+1}(0,q-m+1,\dots,q)\sign R_{m}(\nu).
\end{equation*}
\end{lem}

\begin{proof}[Proof of Theorem~\ref{theo7}]
Let the zeros of $P_{q}(\lambda)$ lie on the interval $(a,b)$ and let $\eta_{b}=P_{q+1}(b)/P_{q}(b)$.
Then the zeros of the polynomial $\widetilde{P}_{q+1}(\lambda)=P_{q+1}(\lambda)-\eta P_{q}(\lambda)$ are such that $\lambda_2,\dots,\lambda_q\in (a,b)$, $\lambda_1<b$ for $\eta<\eta_b$, $\lambda_1=b$ for $\eta=\eta_b$, and $\lambda_1>b$ for $\eta>\eta_b$ (see \cite[Ch.~III, \S\,3.3]{Sz59}).

\smallbreak
{(b)} {In this case $m=0$.} From 
the interlacing of zeros of $P_{\nu}(\lambda)$ (see \eqref{eq3.1}), 
there holds
\[
\lambda_{1,1}<\dots<\lambda_{1,q}<\lambda_{1},
\]
and therefore, $P_{\nu}(\lambda_1)>0$, $\nu\in [0,q]_{\Z}$. Applying the Christoffel-Darboux formula \eqref{eq3.11}, we obtain

\begin{align*}
\sum_{l=0}^{q}d_{l}P_{l}(\lambda_1)P_l(\lambda)&=w_q\,\frac{P_{q+1}(\lambda)P_q(\lambda_1)-P_{q+1}(\lambda_1)P_q(\lambda)}{\lambda-\lambda_1}\\&
=w_q\,\frac{\widetilde{P}_{q+1}(\lambda)P_q(\lambda_1)-\widetilde{P}_{q+1}(\lambda_1)P_{q}(\lambda)}{\lambda-\lambda_1}
=w_qP_q(\lambda_1)\,\frac{\widetilde{P}_{q+1}(\lambda)}{\lambda-\lambda_1}.
\end{align*}
Therefore, the coefficients $a_{\nu}$ in expansion \eqref{eq1.18} satisfy 
\begin{equation}\label{eq4.3}
a_{\nu}=\frac{P_{\nu}(\lambda_1)}{w_qP_q(\lambda_1)}>0,\quad \nu\in [0,q]_{\Z}.
\end{equation}
For $\eta=\eta_b$, we have $\lambda_1=b$ and \eqref{eq1.20} is valid.

{(a, c): the case $m=0$ and  $\eta\neq\eta_b$.} By  the Christoffel-Darboux formula, 
\begin{equation*}
\frac{P_{\nu}(\lambda_1)}{P_{\nu}(b)}-\frac{P_{\nu+1}(\lambda_1)}{P_{\nu+1}(b)}=
\frac{b-\lambda_1}{w_{\nu}P_{\nu}(b)P_{\nu+1}(b)}\sum_{l=0}^{\nu}d_{l}P_{l}(b)P_{l}(\lambda_1).
\end{equation*}
Since $\sign{}(b-\lambda_1)=\sign{}(\eta_b-\eta)$,   we have,
for $\nu\in [0,q-1]_{\Z}$,
\[
\sign \Bigl(\frac{P_{\nu}(\lambda_1)}{P_{\nu}(b)}-\frac{P_{\nu+1}(\lambda_1)}{P_{\nu+1}(b)}\Bigr)=\sign{}(\eta_b-\eta).
\]
From this and  \eqref{eq4.3} we obtain  \eqref{eq1.19} and \eqref{eq1.21} for $m=0$.

\smallbreak
{(a): the case $m\ge 1$ and    $\eta\neq\eta_b$.} Let us calculate the coefficients $a_{\nu}$ in expansion~\eqref{eq1.18} in terms of  $D_{m+1}(\nu,\nu_1,\dots,\nu_m)$, $1 \le\nu_1<\dots<\nu_{m}\le q$.

{Letting $\omega_{m+1}(\lambda)=\prod_{j=1}^{m+1}(\lambda-\lambda_j)$, we have}
\[
\frac{\widetilde{P}_{q+1}(\lambda)}{\omega_{m+1}(\lambda)}
=\sum_{i=1}^{m+1}\frac{\widetilde{P}_{q+1}(\lambda)}{\omega_{m+1}'(\lambda_i)(\lambda-\lambda_i)},
\]
where 
$$\omega_{m+1}'(\lambda_i)=\prod_{j\neq i}(\lambda_i-\lambda_j),\quad \sign \omega_{m+1}'(\lambda_i)=(-1)^{i-1}.$$
Applying the Christoffel-Darboux formula, we obtain as above
\begin{align}
\sum_{\nu=0}^qd_{\nu}P_{\nu}(\lambda_i)P_{\nu}(\lambda)&=w_q\,\frac{P_{q+1}(\lambda)P_q(\lambda_i)-P_{q}(\lambda)P_{q+1}(\lambda_i)}{\lambda-\lambda_i}\notag
\\
&=w_q\,\frac{\widetilde{P}_{q+1}(\lambda)P_q(\lambda_i)}{\lambda-\lambda_i}.\label{eq4.5}
\end{align}
Therefore,
\[
\frac{\widetilde{P}_{q+1}(\lambda)}{\omega_{m+1}(\lambda)}=\frac{1}{w_q}\sum_{\nu=0}^{q-m}d_{\nu}
\sum_{i=1}^{m+1}\frac{P_{\nu}(\lambda_i)}{\omega_{m+1}'(\lambda_i)P_{q}(\lambda_i)}\,P_{\nu}(\lambda)
\]
and
\begin{equation}\label{eq4.6}
a_{\nu}=\frac{1}{w_q}\sum_{i=1}^{m+1}\frac{P_{\nu}(\lambda_i)}{\omega_{m+1}'(\lambda_i)P_{q}(\lambda_i)}.
\end{equation}
{Since}
\[
\Delta(\lambda_1,\dots,\lambda_{m+1})=\prod_{1\le j<s\le m+1}(\lambda_s-\lambda_j)
\]
is  the Vandermonde determinant, {we derive}
\begin{align*}
\frac{1}{\omega_{m+1}'(\lambda_i)}&
=\frac{1}{\Delta(\lambda_1,\dots,\lambda_{m+1})}\,
\frac{\Delta(\lambda_1,\dots,\lambda_{m+1})}{\prod_{j\neq i}(\lambda_i-\lambda_j)}
\\
&=\frac{(-1)^{i-1}\Delta(\lambda_1,\dots,\lambda_{i-1},\lambda_{i+1},\dots,\lambda_{m+1})}{\Delta(\lambda_1,\dots,\lambda_{m+1})}.
\end{align*}
Taking into account  \eqref{eq4.6}, this gives
\begin{align}
a_{\nu}&=\frac{1}{w_q\Delta(\lambda_1,\dots,\lambda_{m+1})}\sum_{i=1}^{m+1}(-1)^{i-1}\Delta(\lambda_1,\dots,\lambda_{i-1},\lambda_{i+1},\dots,\lambda_{m+1})
\frac{P_{\nu}(\lambda_i)}{P_{q}(\lambda_i)}\notag\\ \label{eq4.7}
&=\frac{1}{w_q\Delta(\lambda_1,\dots,\lambda_{m+1})}
\begin{vmatrix}
1&\dots&1\\
\lambda_1&\dots&\lambda_{m+1}\\
\hdotsfor[2]{3}\\
\lambda_1^{m-1}&\dots&\lambda_{m+1}^{m-1}\\
\frac{P_{\nu}(\lambda_1)}{P_{q}(\lambda_1)}&\dots&\frac{P_{\nu}(\lambda_{m+1})}{P_{q}(\lambda_{m+1})}
\end{vmatrix}.
\end{align}

We now rewrite   formula 
\eqref{eq1.10}
as follows:
\[
P_{l-1}(\lambda)=(A_l\lambda+B_l)P_l(\lambda)-C_lP_{l+1}(\lambda),
\]
where
\[
A_l=\frac{\rho_l}{\gamma_{l-1}}>0,\quad B_l=-\frac{\alpha_l}{\gamma_{l-1}},\quad C_l=\frac{\beta_l}{\gamma_{l-1}}>0.
\]
This implies the recurrence relation
\[
\frac{P_{q-j}(\lambda_i)}{P_{q}(\lambda_i)}=(A_{q-j+1}\lambda_i+B_{q-j+1})\,\frac{P_{q-j+1}(\lambda_i)}{P_{q}(\lambda_i)}-
C_{q-j+1}\,\frac{P_{q-j+2}(\lambda_i)}{P_{q}(\lambda_i)}, \quad j\in [1,q]_{\Z},
\]
which yields 
\[
\frac{P_{q-j}(\lambda_i)}{P_{q}(\lambda_i)}=\sum_{s=0}^j\alpha_{s,j}\lambda_i^s,\quad \alpha_{j,j}=A_qA_{q-1}\cdots A_{q-j+1}>0,
\]
where the coefficients $\alpha_{s,j}$ are independent of $\lambda_i$. Substituting these relations into \eqref{eq4.7}, we~get
\begin{align*}
a_{\nu}&=\frac{1}{w_q\Delta(\lambda_1,\dots,\lambda_{m+1})\prod_{j=1}^{m-1}\alpha_{j,j}}
\begin{vmatrix}
\frac{P_{q}(\lambda_1)}{P_{q}(\lambda_1)}&\dots&\frac{P_{q}(\lambda_{m+1})}{P_{q}(\lambda_{m+1})}\\
\frac{P_{q-1}(\lambda_1)}{P_{q}(\lambda_1)}&\dots&\frac{P_{q-1}(\lambda_{m+1})}{P_{q}(\lambda_{m+1})}\\
\hdotsfor[2]{3}\\
\frac{P_{q-m+1}(\lambda_1)}{P_{q}(\lambda_1)}&\dots&\frac{P_{q-m+1}(\lambda_{m+1})}{P_{q}(\lambda_{m+1})}\\
\frac{P_{\nu}(\lambda_1)}{P_{q}(\lambda_1)}&\dots&\frac{P_{\nu}(\lambda_{m+1})}{P_{q}(\lambda_{m+1})}
\end{vmatrix}\\&
=\frac{(P_{q}(\lambda_1)\cdots  P_{q}(\lambda_{m+1}))^{-1}}{w_q\Delta(\lambda_1,\dots,\lambda_{m+1})\prod_{j=1}^{m-1}\alpha_{j,j}}
\begin{vmatrix}
P_{q}(\lambda_1)&\dots&P_{q}(\lambda_{m+1})\\
P_{q-1}(\lambda_1)&\dots&P_{q-1}(\lambda_{m+1})\\
\hdotsfor[2]{3}\\
P_{q-m+1}(\lambda_1)&\dots&P_{q-m+1}(\lambda_{m+1})\\
P_{\nu}(\lambda_1)&\dots&P_{\nu}(\lambda_{m+1})
\end{vmatrix}.
\end{align*}
From the interlacing of zeros
of the polynomials $P_q$ and $\widetilde{P}_{q+1}$, 
we have 
 $\sign P_q(\lambda_i)=(-1)^{i-1}$,
\[
\sign \Delta(\lambda_1,\dots,\lambda_{m+1})=\sign{}(P_{q}(\lambda_1)\cdots P_{q}(\lambda_{m+1}))=(-1)^{\frac{m(m+1)}{2}},
\]
and
\begin{equation}\label{eq4.8}
a_{\nu}=c(q,m)
\begin{vmatrix}
P_{q}(\lambda_1)&\dots&P_{q}(\lambda_{m+1})\\
P_{q-1}(\lambda_1)&\dots&P_{q-1}(\lambda_{m+1})\\
\hdotsfor[2]{3}\\
P_{q-m+1}(\lambda_1)&\dots&P_{q-m+1}(\lambda_{m+1})\\
P_{\nu}(\lambda_1)&\dots&P_{\nu}(\lambda_{m+1})
\end{vmatrix}=c(q,m)D_{m+1}(\nu,q-m+1,\dots,q)
\end{equation}
{with $c(q,m)>0$.}

By Lemma~\ref{lem14},  $D_{m+1}(\nu,q-m+1,\dots,q)$  preserves the sign
for $\nu\in [0,q-m]_{\Z}$
and then, according to \eqref{eq4.8}, all $a_{\nu}$ have the same sign. Further, since the polynomial \eqref{eq1.18} is positive for sufficiently large $\lambda$,
 it follows that
$a_{q-m}>0$ and thus
all $a_{\nu}>0$.
In particular, we have 
$D_{m+1}(0,q-m+1,\dots,q)>0$, which yields, 
by  Lemma~\ref{lem14},
\begin{equation}\label{D-pos}
D_{m+1}(\nu,\nu_1,\dots,\nu_m)>0\quad
\text{for all sequences}\quad 0\le \nu<\nu_1<\dots<\nu_m\le q.
\end{equation}

To show \eqref{eq1.19},
let us write  the difference of coefficients $a_{\nu}$ in \eqref{eq1.18} using \eqref{eq4.8}:
\begin{align*}
J_{\nu}&=\frac{a_{\nu}}{P_{\nu}(b)}-\frac{a_{\nu+1}}{P_{\nu+1}(b)}=c(q,m)
\begin{vmatrix}
P_{q}(\lambda_1)&\dots&P_{q}(\lambda_{m+1})\\
P_{q-1}(\lambda_1)&\dots&P_{q-1}(\lambda_{m+1})\\
\hdotsfor[2]{3}\\
P_{q-m+1}(\lambda_1)&\dots&P_{q-m+1}(\lambda_{m+1})\\
\frac{P_{\nu}(\lambda_1)}{P_{\nu}(b)}-\frac{P_{\nu+1}(\lambda_1)}{P_{\nu+1}(b)}&\dots&
\frac{P_{\nu}(\lambda_{m+1})}{P_{\nu}(b)}-\frac{P_{\nu+1}(\lambda_{m+1})}{P_{\nu+1}(b)}
\end{vmatrix}\\
&=c(q,m)\Bigl(\prod_{j=q-m+1}^{q}P_j(b)\Bigr)\begin{vmatrix}
\frac{P_{q}(\lambda_1)}{P_{q}(b)}&\dots&
\frac{P_{q}(\lambda_{m+1})}{P_{q}(b)}\\
\frac{P_{q-1}(\lambda_1)}{P_{q-1}(b)}-\frac{P_{q}(\lambda_1)}{P_{q}(b)}&\dots&
\frac{P_{q-1}(\lambda_{m+1})}{P_{q-1}(b)}-\frac{P_{q}(\lambda_{m+1})}{P_{q}(b)}\\
\hdotsfor[2]{3}\\
\frac{P_{q-m+1}(\lambda_1)}{P_{q-m+1}(b)}-\frac{P_{q-m+2}(\lambda_1)}{P_{q-m+2}(b)}&\dots&
\frac{P_{q-m+1}(\lambda_{m+1})}{P_{q-m+1}(b)}-\frac{P_{q-m+2}(\lambda_{m+1})}{P_{q-m+2}(b)}\\
\frac{P_{\nu}(\lambda_1)}{P_{\nu}(b)}-\frac{P_{\nu+1}(\lambda_1)}{P_{\nu+1}(b)}&\dots&
\frac{P_{\nu}(\lambda_{m+1})}{P_{\nu}(b)}-\frac{P_{\nu+1}(\lambda_{m+1})}{P_{\nu+1}(b)}
\end{vmatrix}.
\end{align*}


Applying formula \eqref{eq4.5} with $\lambda=b$
and $\widetilde{P}_{q+1}(b)=P_{q}(b)(\eta_{b}-\eta)$, and also
the Christoffel-Darboux formula
with $x=b$ and $y=\lambda_i$, we derive

\begin{align*}
&J_{\nu}\,\frac{(\eta_{b}-\eta)P_{\nu}(b)P_{\nu+1}(b)P_{q}(b)\bigl(\prod_{j=q-m+2}^{q}P_{j}(b)\bigr)w_{\nu}
\prod_{j=q-m+1}^{q}w_j}{c(q,m)\prod_{j=1}^{m+1}(b-\lambda_j)}\\
&=\begin{vmatrix}
\sum_{j=0}^{q}d_jP_{j}(b)P_{q}(\lambda_1)&\dots&\sum_{j=0}^{q}d_jP_{j}(b)P_{j}(\lambda_{m+1})\notag\\
\sum_{j=0}^{q-1}d_jP_{j}(b)P_{j}(\lambda_1)&\dots&\sum_{j=0}^{q-1}d_jP_{j}(b)P_{j}(\lambda_{m+1})\notag\\
\hdotsfor[2]{3}\\
\sum_{j=0}^{q-m+1}d_jP_{j}(b)P_{j}(\lambda_1)&\dots&\sum_{j=0}^{q-m+1}d_jP_{j}(b)P_{j}(\lambda_{m+1})\\
\sum_{j=0}^{\nu}d_jP_{j}(b)P_{j}(\lambda_1)&\dots&\sum_{j=0}^{\nu}d_jP_{j}(b)P_{j}(\lambda_{m+1})\\
\end{vmatrix}\\
&=\Bigl(\prod_{j=q-m+2}^q d_jP_j(b)\Bigr)
\begin{vmatrix}
P_{q}(\lambda_1)&\dots&P_{q}(\lambda_{m+1})\notag\\
P_{q-1}(\lambda_1)&\dots&P_{q-1}(\lambda_{m+1})\notag\\
\hdotsfor[2]{3}\\
P_{q-m+2}(\lambda_1)&\dots&P_{q-m+2}(\lambda_{m+1})\notag\\
\sum_{j=\nu+1}^{q-m+1}d_jP_{j}(b)P_{j}(\lambda_1)&\dots&\sum_{j=\nu+1}^{q-m+1}d_jP_{j}(b)P_{j}(\lambda_{m+1})\notag\\
\sum_{j=0}^{\nu}d_jP_{j}(b)P_{j}(\lambda_1)&\dots&\sum_{j=0}^{\nu}d_jP_{j}(b)P_{j}(\lambda_{m+1})\\
\end{vmatrix}.
\end{align*}
Letting
\[
C_1=\frac{\eta_{b}-\eta}{b-\lambda_1},\quad C_2=\frac{P_{\nu}(b)P_{\nu+1}(b)P_{q}(b)w_{\nu}\prod_{j=q-m+1}^{q}w_j}
{c(q,m)\prod_{j=2}^{m+1}(b-\lambda_j)\prod_{j=q-m+2}^{q}d_j},
\]
for $\nu\in [0,q-m-1]_{\Z}$, we have
\begin{equation}
C_1 C_2 J_{\nu}=
\sum_{\nu_2=\nu+1}^{q-m+1}\sum_{\nu_1=0}^{\nu}d_{\nu_1}d_{\nu_2}P_{\nu_1}(b)P_{\nu_2}(b)D_{m+1}(\nu_1,\nu_2,q-m+2,\dots,q).\label{eq4.9}
\end{equation}
Moreover, for $m=1$, $D_{2}(\nu_1,\nu_2,q-m+2,\dots,q)=D_{2}(\nu_1,\nu_2)$.

We note that 
the right-hand side
of 
 \eqref{eq4.9} {(due to \eqref{D-pos})}
 as well as  the factor $C_2$  are positive.
 We have already mentioned that, for $\eta\neq \eta_b$, the factor $C_1$ is also positive, therefore $J_{\nu}>0$ and \eqref{eq1.19} holds for $m\ge 1$.


{(a): the case $m\ge 1$ and   $\eta=\eta_b$.}
Since, according to \cite[Ch.~III, \S\,3.3]{Sz59}, the function $f(\lambda)=P_{q+1}(\lambda)/P_{q}(\lambda)$ on the interval $(\lambda_{1,q},\infty)$ increases and $f'(b)>0$, by the inverse function~rule,
\[
\lim_{\eta\to \eta_b}\frac{b-\lambda_1(\eta)}{\eta_b-\eta}=\frac{1}{f'(b)}>0.
\]
Thus, $C_1>0$ in
 \eqref{eq4.9} and \eqref{eq1.19} is also valid for $\eta =\eta_b$. 

The proof of Theorem~\ref{theo7} is now complete.
\end{proof}

\section{Proof of Theorem~\ref{theo8}}\label{sec5}



As proven in \cite{Iv21}, polynomials \eqref{eq1.25} and \eqref{eq1.27} are the only extremal polynomials, up to a positive constant, 
in
a version of Yudin's problem, 
where the polynomial coefficients are not restricted to be non-negative.
Thus, 
\cite{Iv21} provides 
the upper bounds for  $B_n(\mu,m)$ in \eqref{eq1.24} and \eqref{eq1.26}. 
To prove the lower bounds for  $B_n(\mu,m)$, it remains to check that polynomials \eqref{eq1.25} and \eqref{eq1.27} belong to  $\Pi_n\cap\Pi_{+} (\{U_l\})$.


(a) In this case $n=2q-m+1$. According to Theorem~\ref{theo7} with $\eta=0$,  we have 
\[
\frac{U_{q+1}(t)}{(t-t_{1})\cdots(t-t_{m+1})}\in \Pi_{+}(\{U_l\}).
\]
Applying the Krein property, we note that 
\[
\frac{U_{q+1}^{2}(t)}{(t-t_{1})\cdots(t-t_{m+1})}\in \Pi_n\cap\Pi_{+}(\{U_l\}).
\]

(b): the case $n=2q-m+2$ and the polynomials $U_{l}^{(1)}(t)$ satisfy the Krein condition. Since (see~\cite{Le83})
\[
U_{l}^{(1)}(t)=\frac{\theta_lU_l(t)+(2-\theta_l)U_{l+1}(t)}{1+t},\quad \theta_l\in (0,2),
\]
applying Theorem~\ref{theo7} for  $U_{q+1}^{(1)}(t)$ with $\eta=0$,
 we derive that the polynomial
\begin{align*}
\frac{(1+t)(U_{q+1}^{(1)}(t))^2}{(t-t_{1}^{(1)})\cdots (t-t_{m+1}^{(1)})}&=\sum_{l=0}^{q-m}d_l^{(1)}b_l(1+t)U_{q+1}^{(1)}(t)U_{l}^{(1)}(t)
\\&=\sum_{l=0}^{q-m}d_l^{(1)}b_l\sum_{k=q+1-l}^{q+1+l}c_{q+1,l,k}(1+t)U_{k}^{(1)}(t)\\&=
\sum_{l=0}^{q-m}d_l^{(1)}b_l\sum_{k=q+1-l}^{q+1+l}c_{q+1,l,k}(\theta_kU_k(t){+}(2-\theta_k)U_{k+1}(t))
\end{align*}
belongs to $\Pi_n\cap\Pi_{+}(\{U_l\})$.

(b): the case $n=2q-m+2$ and $\mu(t)$ be an odd function. Then the measure $d\mu(t)$ is even, $U_l(-1)=(-1)^l$, $\theta_l=1$, and
\begin{equation}\label{eq5.1}
U_{l}^{(1)}(t)=\frac{U_{l}(t)+U_{l+1}(t)}{1+t}.
\end{equation}
Hence, according to the Christoffel-Darboux formula \eqref{eq2.1},
\begin{equation}\label{eq5.2}
U_{l}^{(1)}(t)=w_l^{-1}\sum_{\nu=0}^{l}d_{\nu}(-1)^{l+\nu}U_{\nu}(t).
\end{equation}
Applying \eqref{eq5.1} and \eqref{eq5.2}, we get
\[
\int_{-1}^{1}(U_{l}^{(1)}(t))^2\,(1+t)\,d\mu(t)=
\int_{-1}^{1}w_l^{-1}\sum_{\nu=0}^{l}d_{\nu}(-1)^{l+\nu}U_{\nu}(t)(U_{l}(t)+U_{l+1}(t))\,d\mu(t)=w_l^{-1}.
\]
This,  Theorem~\ref{theo7}, and \eqref{eq5.2} imply that  in the expansion
\[
\frac{U_{q+1}^{(1)}(t)}{(t-t_{1}^{(1)})\cdots (t-t_{m+1}^{(1)})}=
\sum_{l=0}^{q-m}w_lb_{l}U_{l}^{(1)}(t)
\]
the coefficients are monotone, that is, $b_0>b_1>\dots>b_{q-m}>0$. Applying \eqref{eq5.2}, we obtain
\[
\sum_{l=0}^{q-m}w_lb_{l}U_{l}^{(1)}(t)=\sum_{l=0}^{q-m}b_{l}\sum_{\nu=0}^{l}d_{\nu}(-1)^{l+\nu}U_{\nu}(t)
=\sum_{\nu=0}^{q-m}d_{\nu}U_{\nu}(t)\sum_{l=\nu}^{q-m}b_{l}(-1)^{l+\nu}=\sum_{\nu=0}^{q-m}d_{\nu}\delta_{\nu}U_{\nu}(t),
\]
where
\[
\delta_{\nu}=\sum_{l=\nu}^{q-m}b_{l}(-1)^{l+\nu}=b_{\nu}-b_{\nu+1}+\dots+(-1)^{\nu+q-m}b_{q-m}\ge 0.
\]
Thus,
\[
\frac{(1+t)(U_{q+1}^{(1)}(t))^2}{(t-t_{1}^{(1)})\cdots(t-t_{m+1}^{(1)})}=
\sum_{\nu=0}^{q-m}d_{\nu}\delta_{\nu}U_{\nu}(t)(U_{q+1}(t)+U_{q+2}(t))\in \Pi_n\cap\Pi_{+}(\{U_l\}),
\]
completing the proof. \qed


\section{Appendix}
\subsection{Calculation of determinants for Chebyshev polynomials}
In this subsection, we calculate the determinants $D(\nu)$ (see  \eqref{eq4.1}) for the  Jacobi polynomials \cite[Ch.~IV, \S\,4.1]{Sz59} given by
\[
U_l^{(\alpha,\beta)}(t)=\frac{P_{l}^{(-\frac{1}{2}+\alpha,-\frac{1}{2}+\beta)}(t)}
{P_{l}^{(-\frac{1}{2}+\alpha,-\frac{1}{2}+\beta)}(1)},\quad
\alpha,\beta\in\{0,1\}.
\]
It will be convenient to rearrange  $\nu_i$ as follows: 
$1\le\nu_{m}<\dots<\nu_{1}\le q$.

{(i)} 
If $\alpha=\beta=0$, 
 $U_{q+1}^{(0,0)}(t)=P_{q+1}^{(-\frac{1}{2},-\frac{1}{2})}(t)=\cos{}((q+1)\arccos t)$ is the Chebyshev polynomial of the first kind and 
$t_j=\cos\frac{\pi(2j-1)}{2q+2}$, $j=1,\dots,q+1$,
are its  zeros. If $x_i=\frac{\pi\nu_i}{2q+2}$, then
\[
D(\nu)=D(\nu,\nu_m,\dots,\nu_1)=\det{}(\cos{}(2j-1)x_i)_{i,j=1}^{m+1}.
\]
Since $\cos{}(2j-1)x_i=\sum_{s=0}^{j-1}c_{s}\cos^{2s+1}x_i$, $c_{j-1}= 2^{2j-2}$, {the determinant} $D(\nu)$ reduces to the Vandermonde determinant
\begin{align*}
D(\nu)&=\det{}(2^{2j-2}(\cos x_i)^{2j-1})_{i,j=1}^{m+1}=2^{m(m+1)}\prod_{i=1}^{m+1}(\cos
x_i)\det{}((\cos^{2}x_i)^{j-1})_{i,j=1}^{m+1}\\
&=2^{m(m+1)}\prod_{i=1}^{m+1}(\cos x_i)\prod_{1\le k<l\le
m+1}(\cos^2x_l-\cos^2x_k)\\
&=2^{m(m+1)}\prod_{i=1}^{m+1}(\cos x_i)\prod_{1\le k<l\le m}(\cos^2x_l-\cos^2x_k)\,\prod_{1\le k\le m}\Bigl(\cos^2\Bigl(\frac{\pi\nu}{2q+2}\Bigr)-\cos^2\Bigl(\frac{\pi\nu_k}{2q+2}\Bigr)\Bigr).
\end{align*}
Taking into account that $0<\cos x_1<\dots<\cos x_{m+1}<1$, we obtain
\[
\sign D(\nu)=\sign \prod_{k=1}^m \Bigl(\cos\Bigl(\frac{\pi\nu}{2q+2}\Bigr)-\cos\Bigl(\frac{\pi\nu_k}{2q+2}\Bigr)\Bigr)=\sign
\prod_{k=1}^m(\nu_k-\nu).
\]

{(ii)} The  polynomial $U_{q+1}^{(0,1)}(t)$ is given by
\[
U_{q+1}^{(0,1)}(t)=P_{q+1}^{(-\frac{1}{2},\frac{1}{2})}(t)=\frac{\sqrt{2}\cos{}((q+3/2)\arccos t)} {\sqrt{1+t}}
\]
and
$t_j=\cos\frac{\pi(2j-1)}{2q+3}$, $j=1,\dots,q+1$, are its zeros. If $x_i=\frac{\pi(\nu_i+1/2)}{2q+3}$, then as above 
\begin{align*}
D(\nu)&=\det \Bigl(\frac{\sqrt{2}\cos{}(2j-1)x_i}{\sqrt{1+t_j}}\Bigr)_{i,j=1}^{m+1}=\prod_{j=1}^{m+1}\Bigl(\frac{2}{1+t_j}\Bigr)^{1/2}
\det \bigl(\cos{}(2j-1)x_i\bigr)_{i,j=1}^{m+1}\\
&=2^{m(m+1)}\prod_{j=1}^{m+1}\Bigl(\frac{2\cos^2 x_j}{1+t_j}\Bigr)^{1/2}\,\prod_{0\le k<l\le m+1}(\cos^2x_l-\cos^2x_k).
\end{align*}

{(iii)} The polynomial
\[
U_{q+1}^{(1,1)}(t)=P_{q+1}^{(\frac{1}{2},\frac{1}{2})}(t)=\frac{\sin((q+2)\arccos t)} {(q+2)\sqrt{1-t^2}},
\]
is the
 Chebyshev polynomial of the second kind 
with the zeros
$t_j=\cos\frac{\pi j}{q+2}$, $j=1,\dots,q+1$.  If $x_i=\frac{\pi(\nu_i+1)}{q+2}$, then
\[
D(\nu)=\det \Bigl(\frac{\sin jx_i}{(\nu_i+1)\sin\frac{\pi j}{q+2}}\Bigr)_{i,j=1}^{m+1}=
\frac{\det{}(\sin jx_i)_{i,j=1}^{m+1}}{\prod_{j=1}^{m+1}(\nu_j+1)\sin\frac{\pi j}{q+2}}.
\]
Since $\sin jx_i=\sin x_i\,\sum_{s=0}^{j-1}c_{s}\cos^{s}x_i$ with $c_{j-1}=2^{j -1}$, 
we arrive at
\begin{align*}
D(\nu)&=2^{\frac{m(m+1)}{2}}\prod_{j=1}^{m+1}\frac{\sin\frac{\pi(\nu_j+1)}{q+2}}{(\nu_j+1)\sin\frac{\pi
j}{q+2}}\,\det{}(\cos^{j-1}x_i)_{i,j=1}^{m+1}\\
&=2^{\frac{m(m+1)}{2}}\prod_{j=1}^{m+1}\frac{\sin\frac{\pi(\nu_j+1)}{q+2}}{(\nu_j+1)\sin\frac{\pi j}{q+2}}\,
\prod_{0\le k<l\le m+1}(\cos x_l-\cos x_k).
\end{align*}

{(iv)} We have 
\[
U_{q+1}^{(1,0)}(t)=P_{q+1}^{(\frac{1}{2},- \frac{1}{2})}=\frac{\sqrt{2}\sin((q+3/2)\arccos t)} {(2q+3)\sqrt{1-t}},
\]
and $t_j=\cos\frac{\pi j}{q+3/2}$, $j=1,\dots,q+1$, are the zeros of $U_{q+1}^{(1,0)}$.
As in the previous case,
for $x_i=\frac{\pi(\nu_i+1/2)}{q+3/2}$, 
\begin{align*}
D(\nu)&=\det \Bigl(\frac{\sqrt{2}\sin jx_i}{(2\nu_i+1)\sqrt{1-t_j}}\Bigr)_{i,j=1}^{m+1}=
\frac{\det{}(\sin
jx_i)_{i,j=1}^{m+1}}{\prod_{j=1}^{m+1}(2\nu_j+1)\sqrt{(1-t_j)/2}}\\
&=2^{\frac{m(m+1)}{2}}\prod_{j=1}^{m+1}\frac{\sin\frac{\pi(\nu_j+1)}{q+2}}{(2\nu_j+1)\sqrt{(1-t_j)/2}}\prod_{0\le k<l\le m+1}(\cos x_l-\cos x_k).
\end{align*}

\subsection{Corollaries of Theorem \ref{theo7} for trigonometric polynomials}
Finally, writing the Chebyshev polynomials
$U_{q+1}^{(0,0)}$
and $U_{q+1}^{(1,1)}$
 in the trigonometric form, 
the monotonicity property of coefficients in  expansion \eqref{eq1.18} given in Theorem \ref{theo7} implies the following results.


\begin{cor}
Suppose $0\le m\le q$, then
\begin{equation}\label{eq6.1}
\frac{\cos{}(q+1)x}{\bigl(\cos x-\cos\frac{\pi}{2q+2}\bigr)\cdots\bigl(\cos x-\cos\frac{\pi(2m+1)}{2q+2}\bigr)}=
\frac{a_{0,m}}{2}+\sum_{\nu=1}^{q-m}a_{\nu,m}\cos\nu x,
\end{equation}
then
\[
a_{0,m}>a_{1,m}>\dots>a_{q-m,m}>0.
\]
\end{cor}

\begin{cor}
Suppose $1\le m\le q$, then
\[
\frac{\sin(q+1)x}{\bigl(\cos x-\cos\frac{\pi}{q+1}\bigr)\cdots\bigl(\cos x-\cos\frac{\pi m}{q+1}\bigr)}=
\sum_{\nu=1}^{q-m+1}\nu b_{\nu,m}\sin\nu x,
\]
then
\[
b_{1,m}>\dots>b_{q-m+1,m}>0.
\]
\end{cor}

It is worth mentioning that  polynomials 
similar to \eqref{eq6.1}
was considered by Yudin in \cite{Yu02}.



\begin{thebibliography}{99}

\bibitem{Ag00}
R.~P.~Agarwal, \emph{Difference equations and inequalities: theory, methods, and applications}, 2nd ed.,
Marcel Dekker, Inc., New York, 2000. \doi{10.1201/9781420027020}

\bibitem{ABGR05}
R.P.~Agarwal, M.~Bohner, S.R.~Grace, and D.~O'Regan, \emph{Discrete
Oscillation Theory}, Hindawi Publ. Corp., New York, 2005.

\bibitem{At64}
F.~V.~Atkinson, \emph{Discrete and continuous boundary problems}, Academic Press, New York, 1964.

\bibitem{Ba84}
A.~G.~Babenko, \emph{An extremal problem for polynomials}, Math. Notes \textbf{35}(1984), no.~3, 181--186.
\doi{10.1007/BF01139914}

\bibitem{BH20}
P.~B\'erard and B.~Helffer, \emph{Sturm's theorem on zeros of linear
combinations of eigenfunctions}, Expositiones Mathematicae \textbf{38}(2020),
no.~1, 27--50; arXiv:1706.08247v4. \doi{10.1016/j.exmath.2018.10.002}

\bibitem{Bo98}
P.A.~Borodin, \emph{Linearity of metric projections on Chebyshev subspaces in $L^1$ and $C$}, Math. Notes, \textbf{63}(1998), no.~6, 717--723. \doi{10.1007/BF02312764}

\bibitem{CK07}
H.~Cohn and A.~Kumar, \emph{Universally optimal distribution of points on
spheres}, J. Amer. Math. Soc. \textbf{20}(2007), no.~1, 99--148. \doi{10.1090/S0894-0347-06-00546-7}

\bibitem{Ch78}
T.~S.~Chihara, \emph{An Introduction to Orthogonal Polynomials}, Gordon and
Breach Science Publishers, New York--London--Paris, 1978.

\bibitem{DK07}
O.~Do\v{s}l\'{y} and W.~Kratz,
\emph{Oscillation theorems for symplectic difference systems}, J.~Difference Equ. Appl. \textbf{13}(2007), no.~7, 585--605. \doi{10.1080/10236190701264776}

\bibitem{DEH19}
O.~Do\v{s}l\'{y}, J.~Elyseeva, and R.~\v{S}.~Hilscher,
\emph{Symplectic difference systems: oscillation and spectral theory}, Pathw. Math., Springer, Cham, 2019. \doi{10.1007/978-3-030-19373-7}

\bibitem{Du79}
C.~B.~Dunham, \emph{Discrete Chebyshev approximation: alternation and the
Remez algorithm}, Z. Angw. Math. Mech. \textbf{58}(1979), 326--328. \doi{10.1002/zamm.19790590710}

\bibitem{DS08}
V.~K.~Dzyadyk and I.~A. Shevchuk,
\emph{Theory of Uniform Approximation of Functions by Polynomials}, de Gruyter, Berlin, 2008. \doi{10.1515/9783110208245}

\bibitem{El05}
S.~Elaydi, \emph{An introduction to difference equations}, Springer, New York, 2005.

\bibitem{EN04}
A.~Eremenko and D.~Novikov, \emph{Oscillation of Fourier integrals with a spectral gap}, J. Math. Pures Appl. \textbf{83}(2004), no.~3, 313--365.
\doi{10.1016/S0021-7824(03)00064-3}

\bibitem{GK02}
F.~Gantmacher and M.~Krein, \emph{Oscillation Matrices and Kernels and Small
Vibrations of Mechanical Systems}, revised ed., AMS Chelsea Publishing, 2002.

\bibitem{GI00}
D.~V.~Gorbachev and V.~I.~Ivanov, \emph{An extremum problem for polynomials
related to codes and designs}, Math. Notes \textbf{67}(2000), no.~4, 433--438.
\doi{10.1007/BF02676398}

\bibitem{GIT20}
D.~Gorbachev, V.~Ivanov, and S.~Tikhonov, \emph{Uncertainty principles for
eventually constant sign bandlimited functions}, SIAM J. Math. Anal. \textbf{52}
(2020), no.~5, 4751--4782. \doi{10.1137/20M131566X}

\bibitem{GIT24}
D.~Gorbachev, V.~Ivanov, and S.~Tikhonov, \emph{Logan's problem for Jacobi
transforms}, Canad. J. Math. \textbf{76}(2024), no.~3,
4751--4782. \doi{10.4153/S0008414X23000275}

\bibitem{Ha78}
P.~Hartman, \emph{Difference equations: disconjugacy, principal solutions, Green's functions, complete monotonicity}, Trans. Am. Math. Soc. \textbf{246}(1978), 1--30.
\doi{10.1090/S0002-9947-1978-0515528-6}

\bibitem{Iv21}
V.~I.~Ivanov, \emph{Yudin--Hermite extremal problems for polynomials}, {Math. Notes}
\textbf{110}(2021), no.~5, 799--805. \doi{10.4213/mzm13265}

\bibitem{KS66}
S.~Karlin and W.~J.~Studden, \emph{Tchebycheff systems: with applications in analysis
and statistics}, John Wiley~\&~Sons, New York, 1966.

\bibitem{KP01}
W.~G.~Kelley and A.~C.~Peterson, \emph{Difference equations:
an introduction with applications}, 2nd ed., Academic press, San Diego, CA, 2001.

\bibitem{Kr02}
W.~Kratz, \emph{Discrete oscillation}, J.~Difference Equ. Appl. \textbf{9}(2003), no.~1, 135--147. \doi{10.1080/10236100309487540}

\bibitem{La72}
P.-J.~Laurent, Approximation et Optimisation, Hermann, Paris, 1972.

\bibitem{Le83}
V.~I.~Levenshtein, \emph{Boundaries for packings of metric spaces and some
applications}, Problems of Cybernetics \textbf{40}(1983), 43--110. (in
Russian)

\bibitem{LS75}
B.~M.~Levitan and I.~S.~Sargsjan, \emph{Introduction to spectral theory:
selfadjoint ordinary differential operators}, Transl. Math. Monogr. \textbf{39}, AMS, Providence, Rhode Island, 1975.

\bibitem{Lo77}
B.~Logan, \emph{Information in the zero crossings of bandpass signals}, Bell Syst. Tech. J. \textbf{56}(1977), no.~4, 487--510.
\doi{10.1002/j.1538-7305.1977.tb00522.x}

\bibitem{Ma12}
Y.~Mao, \emph{Reconstruction of binary functions and shapes from incomplete frequency information}, IEEE Trans Inf. Theory \textbf{58}(2012), no.~6, 3642--3653.
\doi{10.1109/TIT.2012.2190041}

\bibitem{MP15}
M.~Mitkovski and A.~Poltoratski, \emph{On the determinacy problem for measures}, Invent. Math. \textbf{202}(2015), no.~3, 1241--1267.
\doi{10.1007/s00222-015-0588-6}

\bibitem{MU84}
H.~L.~Montgomery and M.~A. Ulrike, \emph{Biased trigonometric polynomials}. Am. Math. Mon. \textbf{114}(2007), no.~9, 804--809.
\doi{10.1080/00029890.2007.11920472}

\bibitem{MZO24}
S.~Morozov, D.~Zheltkov, and A.~Osinsky,  \emph{Accelerated alternating minimization algorithm for low-rank approximations in the Chebyshev norm}, arXiv:2410.05247.
\doi{10.48550/arXiv.2410.05247}

\bibitem{PK25}
V.~Yu.~Protasov, R.~Kamalov, \emph{How do the lengths of switching intervals influence the stability of a dynamical system?}, Automatica \textbf{171}(2025), 111929; arXiv:2312.10506.
\doi{10.1016/j.automatica.2024.111929}

\bibitem{Si05}
B.~Simon, \emph{Sturm oscillation and comparison theorems}, Amrein, Werner O. (ed.) et al.,
Sturm-Liouville theory. Past and present, 29--43. Birkh\"auser, Basel, 2005. \doi{10.1007/3-7643-7359-8_2}

\bibitem{St20}
S.~Steinerberger, \emph{Quantitative projections in the Sturm oscillation theorem}, J. Math. Pures Appl. \textbf{144}(2020), 1--16.
\doi{10.1016/j.matpur.2020.10.004}

\bibitem{Sz59}
G.~Szeg\"o, \emph{Orthogonal Polynomials}, AMS, New York, 1959.

\bibitem{Te00}
G.~Teschl, \emph{Jacobi Operators and Completely Integrable Nonlinear Lattices}, Mathematical Surveys and Monographs 72, Amer. Math. Soc., Providence, 2000.

\bibitem{Ul06}
A.~Ulanovskii, \emph{The Sturm--Hurwitz theorem and its extensions}, J. Fourier Anal. Appl. \textbf{12}(2006), no.~6, 629--643.
\doi{10.1007/s00041-005-5037-2}

\bibitem{ZMT22}
N.~L.~Zamarashkin, S.~V.~Morozov, and E.~E.~Tyrtyshnikov, \emph{On the best
approximation algorithm by low-rank matrices in Chebyshev's norm}, Comput.
Math. Math. Phys. \textbf{62}(2022), 701--718. \doi{10.1134/S0965542522050141}

\bibitem{Yu97}
V.~A.~Yudin, \emph{Code and design}, Discrete Math. Appl. \textbf{7}(1997), no.~2, 147--155. \doi{10.1515/dma.1997.7.2.147}

\bibitem{Yu02}
V.~A.~Yudin, \emph{Positive values of polynomials}, Math. Notes \textbf{72}(2002), no.~3, 440--443. \doi{10.1023/A:1020520009174}

\bibitem{Yu05}
V.~A.~Yudin, \emph{Distribution of the points of a design on the sphere}, Izv. Math. \textbf{69} (2005), no.~5, 1061--1079.
\doi{10.1070/IM2005v069n05ABEH002288}

\end{thebibliography}
\end{document}